\documentclass[a4paper,12 pt]{article}

\usepackage{amsfonts}
\usepackage{graphicx}
\usepackage{amsmath}
\usepackage{amssymb,subfigure}
\usepackage[english]{babel}
\usepackage{}

\date{}

\newtheorem{thm}{Theorem}[section]

\newtheorem{lemm}[thm]{Lemma}
\newtheorem{remark}[thm]{Remark}

\newtheorem{corl}[thm]{Corollary}
\newtheorem{defin}[thm]{Definition}

 \numberwithin{equation}{section}

\def\R{\mathbb R}

\def\N{\mathbb N}
\def\epsilon{\varepsilon}
\def\ds{\displaystyle}

\newcommand{\be}{\begin{equation}}
\newcommand{\ee}{\end{equation}}
\newcommand{\baa}{\begin{array}}
\newcommand{\eaa}{\end{array}}
\newcommand{\ba}{\begin{eqnarray}}
\newcommand{\ea}{\end{eqnarray}}

\begin{document}
\title{\textbf{Homogenization and influence of fragmentation in a biological invasion model}}
\author{Mohammad El Smaily$^{\hbox{\small{ a}}}$, Fran\c cois Hamel$^{\hbox{\small{ a}}}$ and Lionel Roques$^{\hbox{\small{ b}}}$\\
\\
\footnotesize{$^{\hbox{a }}$Universit\'e Aix-Marseille III, LATP, Facult\'e des Sciences et Techniques}\\
\footnotesize{Avenue Escadrille Normandie-Niemen, F-13397 Marseille Cedex 20, France}\\
\footnotesize{$^{\hbox{b }}$ Unit\'e Biostatistique et  Processus Spatiaux, INRA}\\
\footnotesize{Domaine St Paul - Site Agroparc 84914 Avignon Cedex 9, France}}
\maketitle

\centerline{\it{Dedicated to Professor Masayasu Mimura for his 65th birthday}}
\vskip 0.7cm

\begin{abstract}
In this paper, some properties of the minimal speeds of pulsating Fisher-KPP fronts in periodic environments are established. The limit of the speeds at the homogenization limit is proved rigorously. Near this limit, generically, the fronts move faster when the spatial period is enlarged, but the speeds vary only at the second order. The dependence of the speeds on habitat fragmentation is also analyzed in the case of the patch model.
\end{abstract}

%%%%%%%%%%%%%%%%%%%%%%%%%%%%%%%%%%%%%%%%
%%%%%%%%%%%%%%%%%%%%%%%%%%%%%%%%%%%%%%%%

\section{Introduction and main hypotheses\label{intro}}

In homogeneous environments, the probably most used population dynamics reaction-diffusion model  is the Fisher-KPP model \cite{fi,kpp}. In a one-dimensional space, it corresponds to the following equation
\be\label{kpp}
\frac{\partial u}{\partial t}=D \frac{\partial^2 u}{\partial x^2}+ u\ (\mu-\nu
u), \ t>0, \ x\in \R.
\ee
The unknown $u=u(t,x)$ is the population density at time $t$ and position $x$, and the positive constant coefficients $D$, $\mu$ and $\nu$ respectively correspond to the diffusivity (mobility of the individuals), the intrinsic growth rate and the susceptibility to crowding effects.

A natural extension of this model to heterogeneous environments is the Shigesada-Kawasaki-Teramoto model \cite{skt},
\be\label{eqskt}
\frac{\partial u}{\partial t}=\frac{\partial}{\partial x}\left(a_L(x) \frac{\partial u}{\partial x}\right)+ u\ (\mu_L(x)-\nu_L(x)u), \ t>0, \ x\in \R,
\ee
where the coefficients depend on the space variable $x$ in a L-periodic fashion:

\begin{defin}[$L$-periodicity]\label{L-periodic}
Let $L$ be a positive real number. We say that a function $h:\mathbb{R}\rightarrow\mathbb{R}$ is $L$-periodic if
$$\forall\ x\in\mathbb{R},~h(x+L)=h(x).$$
\end{defin}

In this paper, we are concerned with the general equation:
\begin{equation}\label{eqevo}
\frac{\partial u}{\partial t}=\frac{\partial}{\partial x}\left(a_L(x)\frac{\partial u}{\partial x}\right)+ f_L(x,u),\quad t\in\mathbb{R},~x\in\mathbb{R}.
\end{equation}
The diffusion term $a_L$ satisfies $$a_L(x)=a(x/L),$$ where $a$ is a $C^{2,\delta}(\mathbb{R})$ (with $\delta>0$) 1-periodic function that satisfies
\begin{eqnarray}\label{ca2}
\exists\ 0<\alpha_{1}<\alpha_2,~\forall\ x\in\mathbb{R},~\alpha_1\leq a(x)\leq\alpha_2.
\end{eqnarray}
On other hand, the reaction term satisfies $f_L(x,\cdot)=f(x/L,\cdot)$, where $f:=f(x,s)~:\R\times\R_+\to\R$ is 1-periodic in $x$, of class $C^{1,\delta}$ in $(x,s)$ and $C^2$ in $s$. In this setting, both $a_L$ and $f_L$ are L-periodic in the variable $x$. Furthermore, we assume that:
\begin{eqnarray}\label{cf2}\left\{\begin{array}{ll}
\forall\ x\in\R,\quad f(x,0)=0, \\
\exists\ M\ge 0,\ \forall\ s\ge M,\ \forall\ x\in\R,\quad f(x,s)\le 0, \\
\forall\ x\in\R,\ \ s\mapsto f(x,s)/s\hbox{ is decreasing in }s>0.\end{array}\right.
\end{eqnarray}
Moreover, we set $$\mu(x):=\lim_{s\to 0^+}f(x,s)/s,$$and
$$\mu_L(x):=\lim_{s\to 0^+}f_L(x,s)/s=\mu\left(\frac{x}{L}\right).$$
The growth rate $\mu$ may be positive in some regions (favorable regions) or negative in others (unfavorable regions).

The stationary states $p(x)$ of (\ref{eqevo}) satisfy the equation
\be\label{eqsta}
\frac{\partial}{\partial x}\left(a_L(x)\frac{\partial p}{\partial x}\right)+f_L(x,p)=0,~x\in\mathbb{R}.
\ee
Under general hypotheses including those of this paper, and in any space dimension, it was proved in \cite{bhr1} that a necessary and sufficient condition for the existence of a positive and bounded solution $p$ of (\ref{eqsta}) was the negativity of the principal eigenvalue $\rho_{1,L}$ of the linear operator \be\label{L0}
\mathcal{L}_0: \ \Phi \mapsto-(a_L(x)\Phi')'-\mu_L(x)\Phi,
\ee
with periodicity conditions. In this case, the solution $p$ was also proved to be unique, and therefore L-periodic. Actually, it is easy to see that the map $L\mapsto \rho_{1,L}$ is nonincreasing in $L>0$, and even decreasing as soon as $a$ is not constant (see the proof of Lemma \ref{lemma 1}). Furthermore, $\rho_{1,L}\to-\displaystyle{\int_0^1}\mu(x)dx$ as $L\to 0^+$. In this paper, in addition to the above-mentioned hypotheses, we make the assumption that
\be\label{hypl1}
\displaystyle{\int_0^1}\mu(x)dx>0.
\ee
This assumption then guarantees that
$$\forall\ L>0,\quad \rho_{1,L}<0,$$
whence, for all $L>0$, there exists a unique positive periodic and bounded solution $p_L$ of (\ref{eqsta}). Notice that assumption (\ref{hypl1}) is immediately fulfilled if $\mu(x)$ is positive everywhere.

In this work, we are concerned with the propagation of pulsating traveling fronts which are particular solutions of the reaction-diffusion equation (\ref{eqevo}). Before going further on, we recall the definition of such solutions:

\begin{defin}[Pulsating traveling fronts]\label{pulsating traveling fronts}
A function $u\,=\,u(t,x)$ is called a pulsating traveling front propagating from right to left with an effective speed $c\,\neq\,0,$ if $u$ is a classical solution of:
\begin{eqnarray}\label{front}\left\{\begin{array}{ll}
\displaystyle{\frac{\partial u}{\partial t}=
\frac{\partial}{\partial x}\left(a_L(x)\frac{\partial u}{\partial x}\right)+f_L(x,u),\quad
t\in\mathbb{R},~x\in\mathbb{R},}\\
\displaystyle{\forall\, k\in\mathbb{Z},\; \forall\,(t,x)\,\in\,\mathbb{R}\,\times\,\mathbb{R}},
\quad\displaystyle{u(t+\frac{kL}{c},x)\,=\,u(t,x+kL)}  \hbox{,} \\
0\,\leq\,u(t,x)\leq\,p_L(x),\\
\displaystyle{\lim_{x\,\rightarrow\,-\infty}u(t,x)\,=0\;\hbox{and}\; \lim_{x\rightarrow \,+\infty} u(t,x)-p_L(x)\,=\,0}  \hbox{,}\end{array}\right. \label{defpuls}
\end{eqnarray}
where the above limits hold locally in $t.$
\end{defin}

This definition has been introduced in \cite{sk,skt}. It has also been extended in higher dimensions with $p_L\equiv 1$ in \cite{bh} and \cite{x1}, and with $p_L\not\equiv 1$ in \cite{bhr2}.

Under the above assumptions, it follows from \cite{bhr2} that there exists $c^{*}_L>0$ such that pulsating traveling fronts satisfying (\ref{front}) with a speed of propagation $c$ exist if and only if $c\geq c^{*}_L.$ Moreover, the pulsating fronts (with speeds $c\geq c^{*}_L$) are increasing in time $t.$ Further uniqueness and qualitative properties are proved in \cite{h8,hr2}. The value $\displaystyle{c^{*}_L}$ is called the \emph{minimal speed of propagation}. We refer to \cite{bhn1,bhn2,El Smaily,he8,nad2,nrx,nx,w2} for further existence results and properties of the minimal speeds of KPP pulsating fronts. For existence, uniqueness, stability and further qualitative results for combustion or bistable nonlinearities in the periodic framework, we refer to \cite{clm1,clm2,f,he1,he2,hps,m,na,x1,x2,x3,x6}.

\

In the particular case of the Shigesada {\it{et al}} model (\ref{eqskt}), when $a(x)\equiv 1$, the effects of the spatial distribution of the function $\mu_L$ on the existence and global stability of a positive stationary state $p_L$ of equation (\ref{eqskt}) have been investigated both numerically \cite{rs,sk} and theoretically \cite{bhr1,ccL,rh1}. In particular, as already noticed, enlarging the scale of fragmentation, i.e. increasing $L$, was proved to decrease the value of $\rho_{1,L}$. Biologically, this result means that larger scales have a positive effect on species persistence, for species whose dynamics is modelled by the Shigesada {\it{et al}} model.

The effects of  the spatial distribution of the functions $a_L$ and $\mu_L$ on the minimal speed of propagation $c^*_L$  have not yet been investigated rigorously. This is a difficult problem, since the known variational formula for $c^*_L$ bears on non-self-adjoint operators, and therefore, the methods used to analyze the dependence of $\rho_{1,L}$ on fragmentation cannot be used in this situation. However, in the case of model (\ref{eqskt}), when $a_L\equiv 1$, $\nu_L \equiv 1$ and $\mu_L(x)=\mu(x/L)$, for  a 1-periodic function $\mu$ taking only two values, Kinezaki {\it{et al}} \cite{kkts} numerically observed that $c^*_L$ was an increasing function of the parameter $L$. For sinusoidally varying coefficients, the relationships between $c^*_L$ and $L$ have also been investigated formally by Kinezaki, Kawasaki, Shigesada \cite{kks}. The case of a rapidly oscillating coefficient $a_L(x)$, corresponding to small $L$ values,  and the homogenization limit $L\to 0$, have been discussed in \cite{hps} and \cite{x6} for combustion and bistable nonlinearities $f(u)$.

The first  aim of our work is to analyze rigorously the dependence of the speed of propagation $c^*_L$ with respect to $L$, under the general setting of equation (\ref{eqevo}), for small $L$ values. We determine the limit of the minimal speeds $c^*_L$ as $L\to 0^+$ (the homogenization limit), and we also prove that near the homogenization limit, the species tends to propagate faster when the spatial period of the environment is enlarged. Next, in the case of an environment composed of patches of ``habitat" and ``non-habitat", we  consider the dependence of the minimal speed with respect to habitat fragmentation. We prove that fragmentation decreases the minimal speed.

%%%%%%%%%%%%%%%%%%%%%%%%%%%%%%%%%%%%%%%%
%%%%%%%%%%%%%%%%%%%%%%%%%%%%%%%%%%%%%%%%

\section{Main results\label{results}}

In this section, we describe the main results of this paper. Unless otherwise mentioned, we make the assumptions of Section \ref{intro}. The first theorem  gives the limit of $c_{L}^{*}$ as $L$ goes to $0$.

\begin{thm}\label{thm: c*L--->** as L--->0}
Let $c_{L}^{*}$ be the minimal speed of propagation of pulsating traveling fronts sol\-ving~$(\ref{defpuls})$. Then,
\begin{equation}\label{limit as L tends to 0}
\displaystyle\lim_{L\rightarrow0^{+}}c_{L}^{*}=2\sqrt{<\!a\!>_{H}\ <\!\mu\!>_A},
\end{equation}
where $$<\!\mu\!>_A\ =\ \int_{0}^{1}\mu(x)dx\ \hbox{ and }\ <\!a\!>_{H}\ =\ \left(\int_{0}^{1}(a(x))^{-1}dx\right)^{-1}=\ <\!a^{-1}\!>_A^{-1}$$
denote the arithmetic mean of $\mu$ and the harmonic mean of $a$ over the interval $[0,1].$
\end{thm}

Formula (\ref{limit as L tends to 0}) was derived formally in \cite{skt87} for sinusoidally varying coefficients. Theorem~\ref{thm: c*L--->** as L--->0} then provides a generalization of the formula in \cite{skt87} and a rigorous analysis of the homogenization limit for general diffusion and growth rate profiles.

\begin{remark}\label{rem22}{\rm
The previous theorem gives the limit of $c_{L}^{*}$ as $L\to 0$ when the space dimension is $1.$ Theorem 3.3 of El Smaily \cite{El Smaily} answered this issue in any dimensions~$N$, but under an additional assumption of free divergence of the diffusion field (in the one-dimensional case considered here, this assumption reduces to $da/dx=0$ in $\mathbb{R}$). Lastly, we refer to \cite{clm1,clm2,he1} for other homogenization limits with combustion-type nonlinearities.}
\end{remark}

Our second result describes the behavior of the function
$L\mapsto c^{*}_{L},$ for small $L$ values.

\begin{thm}\label{c*L is increasing}
Let $c_{L}^{*}$ be the minimal speed of propagation of pulsating traveling fronts sol\-ving~$(\ref{defpuls})$. Then, the map $L\mapsto c^*_L$ is of class $C^{\infty}$ in an interval $(0,L_0)$ for some $L_0>0$. Furthermore,
\begin{equation}\label{d/dl c*L-->}
\lim_{L\rightarrow0^{+}}\frac{dc^{*}_{L}}{dL}=0
\ee
and
\be\label{derseconde}
\lim_{L\rightarrow0^{+}}\frac{d^2c^{*}_{L}}{dL^2}=\gamma\ge 0.
\end{equation}
Lastly, $\gamma>0$ if and only if the function
$$\frac{\mu}{<\!\mu\!>_A}+\frac{<\!a\!>_H}{a}$$
is not identically equal to $2$.
\end{thm}

\begin{corl}\label{corollaire} Under the notations of Theorem~$\ref{c*L is increasing}$, it follows that if $a$ is constant and $\mu$ is not constant, or if $\mu$ is constant and $a$ is not constant, then $\gamma>0$ and the speeds $c^*_L$ are increasing with respect to $L$ when $L$ is close to $0$.
\end{corl}

\begin{remark}
{\rm The question of the monotonicity of the map $L\mapsto c^{*}_{L}$ had also been studied under different assumptions in \cite{El Smaily} (see Theorem~5.3). The author answered this question for a reaction-advection-diffusion equation over a periodic domain $\Omega\subseteq\mathbb{R}^{N}$, under an additional assumption on the diffusion coefficient (like in Remark~\ref{rem22}, this assumption would mean again in our present setting that the diffusion coefficient $a(x)$ is constant over $\mathbb{R}$). Our result gives the behavior of the minimal speeds of propagation near the homogenization limit for general diffusion and growth rate coefficients. The condition $\gamma>0$ is generically fulfilled, which means that, roughly speaking, the more oscillating the medium is, the slower the species moves. But the speeds vary only at the second order with respect to the period $L$. Based on numerical observations which have been carried out in \cite{kks} for special types of diffusion and growth rate coefficients, we conjecture that the monotonicity of $c^*_L$ holds for all $L>0$.}
\end{remark}

Lastly, we give a first theoretical evidence that habitat fragmentation, without changing the scale $L$, can decrease the minimal speed $c^*$. We here fix a period $L_0>0$.

\begin{figure} \centering
\subfigure[]{
\includegraphics*[width=7cm]{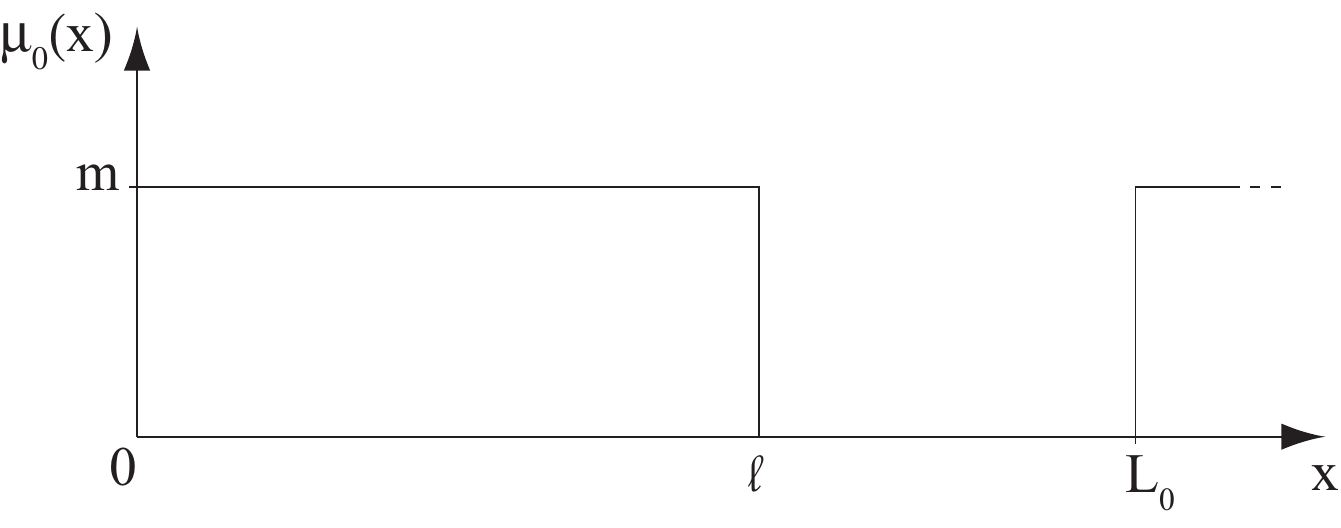}}
\subfigure[]{\includegraphics*[width=7cm]{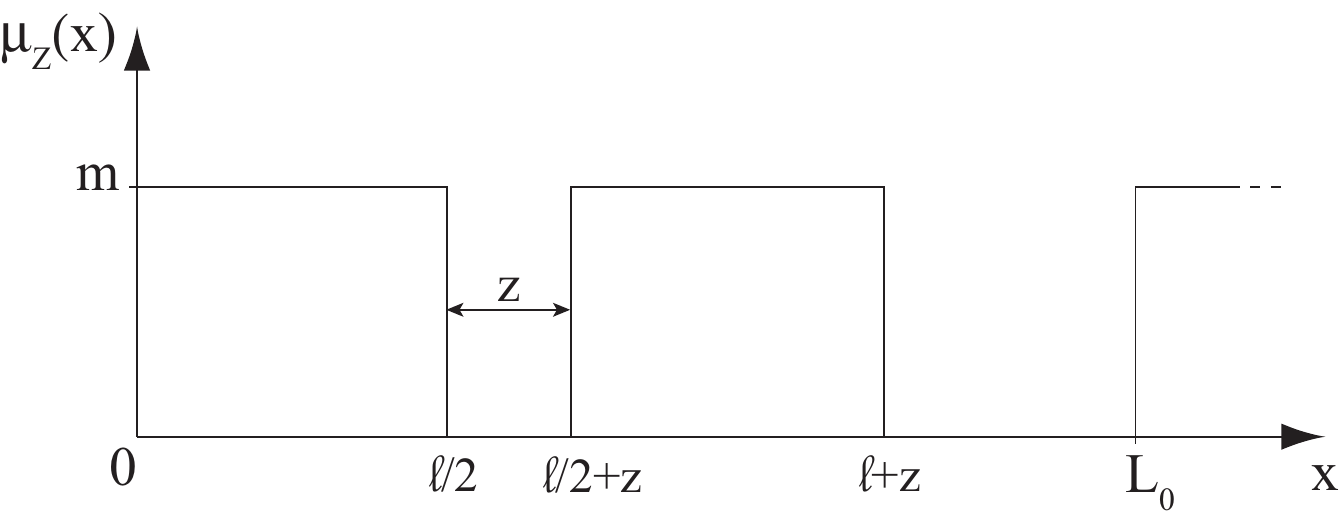}} \caption{The
$L_0$-periodic function $x\mapsto \mu_{z}(x)$, (a): with $z=0$;
(b): with $z>0$.} \label{fig:1}
\end{figure}

We assume that $a\equiv 1$, and that $\mu_{L_0}:=\mu_z$ takes only the two values $0$ and $m>0$, and depends on a parameter $z$. More precisely:
\be\label{defmuz}\left\{\baa{l}
\hbox{There exist }0\leq z \hbox{ and }l\in(0,L_0) \hbox{ such that
}l+z\leq L_0, \\ \mu_z \equiv m  \hbox{ on } [0,l/2)\cup [l/2+z, l+z), \\
\mu_z \equiv 0  \hbox{ on } [l/2,l/2+z)\cup [l+z,L_0).\eaa
\right.\ee
With this setting, the region where $\mu_z$ is positive, which can be interpreted as ``habitat" in the Shigesada {\it{et al}} model, is of Lebesgue measure $l$ in each period cell $[0,L_0]$. For $z=0$, this region is simply an interval. However, whenever $z$ is positive, this region is fragmented into two parts of same length $l/2$ (see Figure \ref{fig:1}). Our next result means that this fragmentation into two parts reduces the speed $c^*$.

\begin{thm} Let $c^*_z$ be the minimal speed of propagation of pulsating traveling fronts sol\-ving~$(\ref{defpuls})$, with $a_{L_0}\equiv 1$ and $\mu_{L_0} = \mu_z$ defined by $(\ref{defmuz})$. Assume that $l\in(3 L_0/4 ,L_0)$. Then  $z\mapsto c^*_z $ is decreasing in $[0,(L_0-l)/2]$, and increasing in $[(L_0-l)/2,L_0-l].$
\label{th_fragmentation}
\end{thm}

\begin{remark}\label{rem31}{\rm Note that, whenever $z>(L_0-l)/2$, the two habitat components in the period cell $[l/2+z,L_0+l/2+z]$ are at a distance smaller than $(L_0-l)/2$ from each other.  In fact, Theorem~\ref{th_fragmentation} proves that, when $z$ varies in $(0,L_0-l)$, $c^*_z $ is all the larger as the minimal distance separating two habitat components is small, that is as the maximal distance between two consecutive habitat components is large.}
\end{remark}

\begin{remark}\label{rem32}{\rm Here, the function $\mu_z$ does not satisfy the general regularity assumptions of Section~\ref{intro}. However, $c^*_z $ can still be interpreted as the minimal speed of propagation of weak solutions of (\ref{defpuls}), whose existence can be obtained by approaching $\mu_z$ with regular functions.}
\end{remark}

The main tool of this paper is a variational formulation for $c^{*}_L$ involving elliptic eigenvalue problems which depend strongly on the coefficients $a$ and $f.$ Such a formulation was given in any space dimension in \cite{bhn2} in the case where the bounded stationary state $p$ of the equation (\ref{eqevo}) is constant, and in \cite{bhr2} in the case of a general nonconstant  bounded stationary state~$p(x)$.

%%%%%%%%%%%%%%%%%%%%%%%%%%%%%%%%%%%%%%%%
%%%%%%%%%%%%%%%%%%%%%%%%%%%%%%%%%%%%%%%%

\section{The homogenization limit: proof of Theorem \ref{thm: c*L--->** as L--->0} \label{proof1}}

This proof is divided into three main steps.

\underline{Step 1: a rough upper bound for $c^*_L$.} For each $L>0,$ the minimal speed $c^{*}_L$ is positive and, from \cite{bhr2} (see also \cite{bhn2}  in the case when $p\equiv 1$), it is given by the variational formula
\begin{equation}\label{var formula}
c^{*}_L=\min_{\lambda>0}\frac{k(\lambda,L)}{\lambda}=\frac{k(\lambda^{*}_L,L)}{\lambda^{*}_L},
\end{equation}
where $\lambda^*_L>0$ and, for each $\lambda\in\R$ and $L>0$, $k(\lambda,L)$ denotes the principal eigenvalue of the problem
\begin{equation}\label{elliptic equation L periodic}
\left(a_{L}\psi_{\lambda,L}'\right)'+2\lambda a_L\psi'_{\lambda,L}+\lambda
a_{L}'\psi_{\lambda,L}+\lambda^{2}a_L\psi_{\lambda,L}+\mu_L\psi_{\lambda,L}=k(\lambda,L)\psi_{\lambda,L}~\hbox{in}~\R,
\end{equation}
with $L$-periodicity conditions.  In (\ref{elliptic equation L periodic}), $\psi_{\lambda,L}$ denotes a principal eigenfunction, which is of class $C^{2,\delta}(\R)$, positive, unique up to multiplication by a positive constant, and $L$-periodic. Furthermore, it follows from Section~3 of \cite{bhr2}
that the map $\lambda\mapsto k(\lambda,L)$ is convex and that $\frac{\partial k}{\partial\lambda}(0,L)=0$ for each $L>0$. Therefore, for each $L>0$, the map $\lambda\mapsto k(\lambda,L)$ is nondecreasing in $\R_+$ and
\be\label{klambdaL}
\forall\ \lambda\ge0,\ \forall\ L>0,\quad k(\lambda,L)\ge k(0,L)=-\rho_{1,L}>0
\ee
under the notations of Section~\ref{intro}.\par
Multiplying (\ref{elliptic equation L periodic}) by $\psi_{\lambda,L}$ and integrating by parts over $[0,L],$ we get, due to the $L$-periodicity of $a_L$ and $\psi_{\lambda,L}$~:
$$k(\lambda,L)\int_0^L\psi_{\lambda,L}^2=-\int_{0}^{L}a_L\left(\psi_{\lambda,L}'\right)^2\;+\lambda^{2}\int_{0}^{L}a_{L}\psi_{\lambda,L}^{2}+\int_{0}^{L}\mu_{L}\psi_{\lambda,L}^{2},$$
for all $\lambda>0$ and for all $L>0.$ Consequently,
\begin{equation}\label{K/lambda}
\forall\ \lambda>0,~\forall\ L>0,\quad k(\lambda,L)\leq\lambda^2a_M+\mu_M,
\end{equation}
where
$$a_M=\max_{x\in\mathbb{R}}a(x)>0\hbox{ and }\mu_M=\max_{x\in\mathbb{R}}\mu(x)>0.$$\par
Using (\ref{var formula}), we get that
\begin{equation}{\label{c*L is bounded}}
\forall\ L>0,~~0<c^{*}_{L}\leq 2\sqrt{a_M\mu_M}.
\end{equation}

\underline{Step 2: the sharp upper bound for $c^*_L$.} For any $\lambda>0$ and $L>0$, consider the functions
$$\varphi_{\lambda,L}(x):=e^{\lambda\,x}\psi_{\lambda,L}(x),~~x\in\mathbb{R}.$$
Since $\psi_{\lambda,L}$ is unique up to multiplication, we will assume in this step~2 that
\begin{equation}\label{varphi in L2}
\int_0^2\varphi_{\lambda,L}^2(x)dx=1.
\ee
The above choice ensures that
\begin{equation}\label{psi L2}
\int_{0}^2\psi_{\lambda,L}^{2}(x)dx\leq\int_{0}^{2}e^{2\lambda x}\psi_{\lambda,L}^{2}(x)dx=\int_{0}^{2}\varphi_{\lambda,L}^{2}(x)dx=1.
\end{equation}\par
We are now going to prove that the families $(\psi_{\lambda,L})_{\lambda,L}$ and $(\varphi_{\lambda,L})_{\lambda,L}$ remain bounded in $H^{1}(0,1)$ for $L$ small enough and as soon as $\lambda$ stays bounded. For each $L>0,$ we call
$$M_L=[1/L]+1\in\mathbb{N},$$
where $[1/L]$ stands for the integer part of $1/L$. Multiplying (\ref{elliptic equation L periodic}) by $\psi_{\lambda,L}$ and integrating by parts over $[0,M_L\,L],$ we get that
$$-\int_{0}^{M_LL}a_L{\psi_{\lambda,L}'}^{2}+\int_{0}^{M_L\,L}\lambda^{2}a_{L}\psi_{\lambda,L}^{2}+\int_{0}^{M_L\,L}\mu_L\psi_{\lambda,L}^{2}=k(\lambda,L)\int_{0}^{M_LL}\psi_{\lambda,L}^{2}.$$
Using (\ref{ca2}), (\ref{klambdaL}) and (\ref{K/lambda}), it follows that
$$\displaystyle{0\le\int_{0}^{M_LL}{\psi_{\lambda,L}'}^{2}\leq\frac{1}{\alpha_1}
\times\left(\lambda^{2}a_M+\mu_M\right)\times\int_{0}^{M_L\,L}
\psi_{\lambda,L}^{2}}.$$
Since $1<M_LL\le 1+L$ for all $L>0$, we have that $1<M_LL\le2$ for all $L\leq 1$. Thus, for all $0<L\leq 1$,
$$\int_{0}^{1}{\psi_{\lambda,L}'}^{2}\leq\int_{0}^{M_LL}{\psi_{\lambda,L}'}^{2}\hbox{ and }\int_0^{M_LL}\psi_{\lambda,L}^2\le\int_0^2\psi_{\lambda,L}^2\le 1$$
from (\ref{psi L2}). It follows now that
\begin{equation}\label{psi' in L2}
\forall\ \lambda>0,\ \forall\ 0<L\leq 1,~~\int_{0}^{1}{\psi_{\lambda,L}'}^{2}\leq\frac{\lambda^{2}a_M+\mu_M}{\alpha_1}.
\end{equation}
From (\ref{psi L2}) and (\ref{psi' in L2}), we conclude that, for any given $\Lambda>0$, the family $\displaystyle{(\psi_{\lambda,L})_{0<\lambda\le\Lambda,\ 0<L\leq 1}}$ is bounded in $H^{1}(0,1).$ On the other hand,
$$\varphi_{\lambda,L}'(x)=\lambda\varphi_{\lambda,L}(x)+e^{\lambda
x}\psi_{\lambda,L}'(x).$$
Owing to (\ref{varphi in L2}) and (\ref{psi' in L2}), we get:
\begin{eqnarray}\label{varphi' L2}
%\nonumber to remove numbering (before each equation)
\begin{array}{ll}
\displaystyle{\forall\ \lambda>0,\ \forall\ L\leq 1,~~||\varphi_{\lambda,L}'||_{L^{2}(0,1)}}&\displaystyle{\leq\lambda\underbrace{||\varphi_{\lambda,L}||_{L^{2}(0,1)}}_{\leq1}+e^{\lambda}||\psi_{\lambda,L}'||_{L^{2}(0,1)}}\\
&\displaystyle{\leq\lambda+e^{\lambda}\times\sqrt{\frac{\lambda^{2}a_M+\mu_M}{\alpha_1}}}.
\end{array}
\end{eqnarray}
From (\ref{varphi in L2}) and (\ref{varphi' L2}), we obtain that, for any given $\Lambda>0$, the family $\displaystyle{(\varphi_{\lambda,L})_{0<\lambda\le\Lambda,\ 0<L\leq 1}}$ is bounded in $H^{1}(0,1)$ and that the family $\displaystyle{(a_L\varphi'_{\lambda,L})_{0<\lambda\le\Lambda,\ 0<L\leq 1}}$ is bounded in  $L^{2}(0,1)$ (due to~(\ref{ca2})). Moreover,
$$\left(a_L\varphi_{\lambda,L}'\right)'=\lambda^{2}a_L e^{\lambda
\,x}\psi_{\lambda,L}+2\lambda a_Le^{\lambda \,x}\psi'_{L}+\lambda
a_{L}'e^{\lambda \,x}\psi_{\lambda,L}+ e^{\lambda
\,x}a_{L}'\psi_{\lambda,L}'+e^{\lambda \,x}a_{L}\psi_{\lambda,L}''.$$
Multiplying (\ref{elliptic equation L periodic}) by $e^{\lambda\,x},$ we then get
\begin{equation}\label{multiply ell by exp x}
\left(a_L\varphi_{\lambda,L}'\right)'+\mu_{L}\varphi_{\lambda,L}=k(\lambda,L)\varphi_{\lambda,L}~\hbox{in}~\mathbb{R}.
\end{equation}
Let
$$v_{\lambda,L}(x)=a_L(x)\varphi_{\lambda,L}'(x)$$
for all $\lambda>0$, $L>0$ and $x\in\mathbb{R}.$ Pick any $\Lambda>0$. One already knows that the family $(v_{\lambda,L})_{0<\lambda\le\Lambda,\ 0<L\leq 1}$ is bounded in $L^{2}(0,1).$ Furthermore,
\begin{equation}\label{equation sat by vL}
v_{\lambda,L}'+\mu_{L}\varphi_{\lambda,L}=k(\lambda,L)\varphi_{\lambda,L}~\hbox{in}~\mathbb{R}.
\end{equation}
Notice that the family $(k(\lambda,L))_{0<\lambda\le\Lambda,\ 0<L\le 1}$ is bounded from (\ref{klambdaL}) and (\ref{K/lambda}). From (\ref{varphi in L2}) and  (\ref{equation sat by vL}), it follows that the family $\left(v_{\lambda,L}'\right)_{0<\lambda\le\Lambda,\ 0<L\leq 1}$ is bounded in $L^2(0,1).$ Eventually, $(v_{\lambda,L})_{0<\lambda\le\Lambda,\ 0<L\leq 1}$ is bounded in $H^{1}(0,1)$.\par
Pick now any sequence $(L_n)_{n\in\N}$ such that $0<L_n\le 1$ for all $n\in\N$, and $L_n\rightarrow0^+$ as $n\rightarrow+\infty.$ Choose any $\lambda>0$ and any sequence $(\lambda_n)_{n\in\N}$ of positive numbers such that $\lambda_n\to\lambda$ as $n\to+\infty$. We claim that
\be\label{claim}
k(\lambda_n,L_n)\to\lambda^{2}<\! a\!>_H+<\!\mu\!>_A\hbox{ as }n\to+\infty,
\ee
where $<\! a\!>_H=\displaystyle{\left(\int_0^1(a(x))^{-1}dx\right)}$ and $<\!\mu\!>_A=\displaystyle{\int_0^1}\mu(x)dx$. To do so, call
$$\psi_n=\psi_{\lambda_n,L_n},\ \varphi_n=\varphi_{\lambda_n,L_n}\hbox{ and }v_n=v_{\lambda_n,L_n}.$$
It follows from the above computations that the sequences $(\psi_n)$ and $(v_n)$ are bounded in $H^{1}(0,1).$ Hence, up to extraction of a subsequence,
$$\psi_n\rightarrow\overline{\psi}\hbox{ and }v_n\rightarrow w\hbox{ as }n\rightarrow+\infty,$$
strongly in $L^2(0,1)$ and weakly in $H^{1}(0,1)$. By Sobolev injections, the sequence $(\psi_n)$ is bounded in $C^{0,1/2}([0,1])$. But since each function $\psi_n$ is $L_n$-periodic (with $L_n\to 0^+$), it follows from Arzela-Ascoli theorem that $\overline{\psi}$ has to be constant over $[0,1]$. Moreover, the boundedness of the sequence $\left(k(\lambda_n,L_n)\right)_{n\in\N}$ implies that, up to extraction of another subsequence,
$$k(\lambda_n,L_n)\to\overline{k}(\lambda)\in\R\hbox{ as }n\rightarrow+\infty.$$
We denote this limit by $\overline{k}(\lambda)$, we will see later that indeed it depends only on $\lambda$. It follows now, from (\ref{equation sat by vL}) after replacing $(\lambda,L)$ by $(\lambda_n,L_n)$ and passing to the limit as $n\rightarrow+\infty$, that
$$w'+<\!\mu\!>_Ae^{\lambda x}\overline{\psi}=\overline{k}(\lambda)\overline{\psi}e^{\lambda x}~~\hbox{a.e. in}~(0,1).$$
Notice indeed that $\mu_L\rightharpoonup<\!\mu\!>_A$ as $L\to0^+$ in $L^2(0,1)$ weakly. Meanwhile,
$$\varphi'_{n}=\lambda_ne^{\lambda_nx}\psi_{n}+e^{\lambda_nx}{\psi'_{n}}=
\frac{v_{n}}{a_{L_n}}\rightharpoonup\ <\! a^{-1}\!>_A\ w\hbox{ as }n\to+\infty,\hbox{ weakly in }L^2(0,1),$$
where $<\! a^{-1}\!>_A=\displaystyle{\int_{0}^{1}}(a(x))^{-1}dx$. Thus, we obtain
$$w=<\! a^{-1}\!>_A^{-1}\lambda e^{\lambda x}\overline{\psi}=<\! a\!>_H\lambda e^{\lambda x}\overline{\psi}.$$
Consequently,
$$\displaystyle{\lambda^{2}<\! a\!>_H\,\overline{\psi}}+<\!\mu\!>_A\;\overline{\psi}=\overline{k}(\lambda)\overline{\psi}.$$
Actually, since the functions $\psi_n$ are $L_n$-periodic (with $L_n\to 0^+$) and converge to the constant $\overline{\psi}$ strongly in $L^2(0,1)$, they  converge to $\overline{\psi}$ in $L^2_{loc}(\R)$. But
$$1=\int_0^2\varphi_n^2\le e^{4\lambda_n}\int_0^2\psi_n^2\le e^{4M}\int_0^2\psi_n^2,$$
where $M=\sup_{n\in\N}\lambda_n$. Hence, $\overline{\psi}\neq 0$ and
\begin{equation}\label{msh}
\lambda^{2}<\! a\!>_H+<\!\mu\!>_A=\overline{k}(\lambda).
\end{equation}
By uniqueness of the limit, one deduces that the whole sequence $(k(\lambda_n,L_n))_{n\in\N}$ converges to this quantity $\overline{k}(\lambda)$ as $n\to+\infty$, which proves the claim (\ref{claim}).\par
Now, take any sequence $L_n\to 0^+$ such that $c^*_{L_n}\to\limsup_{L\to 0^+}c^*_L$ as $n\to+\infty$. For each $\lambda>0$ and for each $n\in\N$, one has
$$c^*_{L_n}\le\frac{k(\lambda,L_n)}{\lambda}$$
from (\ref{var formula}), whence
$$\limsup_{L\to 0^+}c^*_L=\lim_{n\to+\infty}c^*_{L_n}\le\frac{\overline{k}(\lambda)}{\lambda}=\lambda<\! a\!>_H+\frac{<\!\mu\!>_A}{\lambda}.$$
Since this holds for all $\lambda>0$, one concludes that
\begin{equation}\label{lim sup}
\limsup_{L\rightarrow0^{+}}c^{*}_{L}\leq 2\sqrt{<\! a\!>_H\ <\!\mu\!>_A}.
\end{equation}

\underline{Step 3: the sharp lower bound for $c^*_L$.} The aim of this step is to prove that
$$\liminf_{L\rightarrow0^{+}}c^{*}_{L}\geq 2\sqrt{<\! a\!>_H\ <\!\mu\!>_A}$$
which would complete the proof of Theorem~\ref{thm: c*L--->** as L--->0}.\par
For each $L>0,$ the minimal speed $c^{*}_{L}$ is given by (\ref{var formula}) and the map $(0,+\infty)\ni\lambda\mapsto k(\lambda,L)/\lambda$ attains its minimum at $\lambda^*_L>0$.We will prove that, for $L$ small enough, the family $(\lambda^{*}_L)$ is bounded from above and from below by $\overline{\lambda}>0$ and $\underline{\lambda}>0$ respectively. Namely, one has

\begin{lemm}\label{lemma 1} There exist $L_0$ and $0<\underline{\lambda}\le\overline{\lambda}<+\infty$ such that
$$\underline{\lambda}\;\leq\lambda^{*}_L\;\leq\overline{\lambda}~~\hbox{for all}~~0<L\leq L_0.$$
\end{lemm}

The proof is postponed at the end of this section. Take now any sequence $(L_n)_n$ such that $0<L_n\le L_0$ for all $n,$ and $L_n\rightarrow0^+$ as $n\rightarrow+\infty.$ From Lemma~\ref{lemma 1}, there exists $\lambda^{*}>0$ such that, up to extraction of a subsequence, $\lambda^{*}_{L_n}\rightarrow\lambda^{*}$ as $n\rightarrow+\infty.$ One also has
$$c^{*}_{L_n}=\displaystyle{\frac{k(\lambda^{*}_{L_n},L_n)}{\lambda^{*}_{L_n}}}
\displaystyle{\mathop{\rightarrow}_{n\to+\infty}}\displaystyle{\frac{\overline{k}(\lambda^{*})}{\lambda^{*}}=\lambda^*<\! a\!>_H+\frac{<\!\mu\!>_A}{\lambda^*}\ge2\sqrt{<\! a\!>_H\ <\!\mu\!>_A}}$$
from (\ref{claim}) and (\ref{msh}). Therefore, $\liminf_{L\to 0^+}c^*_L\ge2\sqrt{<\! a\!>_H\ <\!\mu\!>_A}$. Eventually,
$$\lim_{L\rightarrow0^{+}}c^{*}_{L}=2\sqrt{<\! a\!>_H\ <\!\mu\!>_A}$$
and the proof of Theorem~\ref{thm: c*L--->** as L--->0} is complete.\hfill$\Box$\break

\noindent\textbf{Proof of Lemma~\ref{lemma 1}.} Observe first that, for $\lambda=0$ and for any $L>0$, $k(0,L)$ is the principal eigenvalue of the problem
$$(a_L\phi_L')'+\mu_L\phi_L=k(0,L)\phi_L~~\hbox{in}~\mathbb{R},$$
and we denote $\phi_L=\psi_{0,L}$ a principal eigenfunction, which is $L$-periodic, positive and unique up to multiplication. In other words, $k(0,L)=-\rho_{1,L}$ under the notations of Section~\ref{intro}. Dividing the above elliptic equation by $\phi_L$ and integrating by parts over $[0,L],$ one gets
$$k(0,L)=\frac{1}{L}\int_{0}^{L}\frac{a_L\,{\phi'_L}^{2}}{\phi_L^2}+\int_{0}^{1}\mu(x)dx\ \ge\ <\!\mu\!>_A\ >0.$$
On the other hand, as already recalled, $\frac{\partial k}{\partial\lambda}(0,L)=0$  and the map $\lambda\mapsto k(\lambda,L)$ is convex for all $L>0$. Therefore,
$$\forall\ \lambda>0,~\forall\ L>0,~k(\lambda,L)\ \geq\ k(0,L)\ \geq\ <\!\mu\!>_A\ >\ 0.$$\par
Assume here that there exists a sequence $(L_n)_{n\in\N}$ of positive numbers such that $L_n\to 0^+$ and $\lambda^{*}_{L_n}\rightarrow0^+$ as $n\to+\infty$. One then gets $$c^*_{L_n}=\frac{k(\lambda^{*}_{L_n},L_n)}{\lambda^{*}_{L_n}}\geq\frac{<\!\mu\!>_A}{\lambda^{*}_{L_n}}\to+\infty\hbox{ as }n\to+\infty.$$
This is contradiction with (\ref{lim sup}). Thus, for $L>0$ small enough, the family $(\lambda^{*}_{L})_L$ is bounded from below by a positive constant $\underline{\lambda}>0$ (actually, these arguments show that the whole family $(\lambda^*_L)_{L>0}$ is bounded from below by a positive constant).\par
It remains now to prove that $(\lambda^{*}_L)_{L}$ is bounded from above when $L$ is small enough. We assume, to the contrary, that there exists a sequence $L_n\rightarrow0^+$ as $n\rightarrow+\infty$ such that  $\lambda^{*}_{L_n}\rightarrow+\infty$ as $n\rightarrow+\infty.$ Call
$$k_n=k(\lambda^*_{L_n},L_n),\ \psi_n(x)=\psi_{\lambda^*_{L_n},L_n}(x)\hbox{ and } \varphi_n(x)
=\varphi_{\lambda^*_{L_n},L_n}(x)=e^{\lambda^*_{L_n}x}\psi_n(x)$$
for all $n\in\N$ and $x\in\R$. Rewriting (\ref{multiply ell by exp x}) for $\lambda=\lambda^{*}_{L_n}$ and for $L=L_n,$ one consequently gets
\begin{equation}\label{upper bdd 1}
\forall\ n\in\N,~~\left(a_{L_n}\varphi_{n}'\right)'+\mu_{L_n}\varphi_{n}=
k_n\varphi_{n}~\hbox{in}~\mathbb{R}.
\end{equation}
Owing to the positivity and the $L_n$-periodicity of the $C^{2}(\mathbb{R})$ eigenfunction $\psi_{n},$ it follows that
$$ \forall\ n\in\mathbb{N},~~\exists\ \theta_{n}\in[0,L_n],~~ \psi_{n}(\theta_{n})=\max_{x\in\mathbb{R}}\psi_{n}(x)=\max_{x\in[0,L_n]}\psi_{n}(x),$$
whence
$$\forall n\in\mathbb{N},~~\psi_{n}'(\theta_{n})=0.$$
For each $n\in\N$, let $M_{L_n}=[1/L_n]+1\in\mathbb{N}.$ Thus,
$$\forall\  n\in\mathbb{N},~\varphi_{n}'(\theta_{n}+M_{L_n}L_n)=\lambda^{*}_{L_n}e^{\lambda^{*}_{L_n}(\theta_n+M_{L_n}L_n)}\psi_{n}(\theta_{n}).$$
Multiplying (\ref{upper bdd 1}) by $\varphi_{n}$ and integrating by parts over the interval $[\theta_n,\theta_n+M_{L_n}\,L_n],$ one then obtains
\begin{equation}\label{int by parts over theta_L}\begin{array}{lc}
&\displaystyle{\underbrace{a_{L_n}(\theta_n+M_{L_n}\,L_n)\varphi'_{n}(\theta_n+M_{L_n}\,L_n)\varphi_{n}(\theta_n+M_{L_n}\,L_n)-a_{L_n}(\theta_n)\varphi'_{n}(\theta_n)\varphi_{n}(\theta_n)}_{A(n)}}\\
&\displaystyle{-\underbrace{\int_{\theta_n}^{\theta_n+M_{L_n}\,L_n}a_{L_n}{\varphi'_{n}}^{2}}_{B(n)}+\underbrace{\int_{\theta_n}^{\theta_n+M_{L_n}\,L_n}\mu_{L_n}\varphi_{n}^{2}}_{C(n)}=k_n\int_{\theta_n}^{\theta_n+M_{L_n}\,L_n}\varphi_{n}^{2}}.
\end{array}
\end{equation}
But, for each $n\in\N$, $M_{L_n}\in\mathbb{N}$ while $a_{L_n}$ and $\psi_{n}$ are $L_n$-periodic. Hence, $a_{L_n}(\theta_n+M_{L_n}\,L_n)=a_{L_n}(\theta_n)$, $\psi_{n}(\theta_n+M_{L_n}\,L_n)=\psi_{n}(\theta_n),$ and $\psi'_{n}(\theta_n+M_{L_n}\,L_n)=\psi'_{n}(\theta_n)=0.$ Then,
\begin{equation}\label{A(L)}\begin{array}{lll}
A(n)&=&\displaystyle{a_{L_n}(\theta_n)\lambda^{*}_{L_n}\psi^{2}_{n}(\theta_n)\left(e^{2\lambda^{*}_{L_n}(\theta_n+M_{L_n}\,L_n)}-e^{2\lambda^{*}_{L_n}\theta_n}\right)}\\
&\geq&\displaystyle{\frac{\alpha_{1}}{2}\times\lambda^{*}_{L_n}\psi^{2}_{n}(\theta_n)e^{2\lambda^{*}_{L_n}(\theta_n+M_{L_n}\,L_n)}}~~\hbox{($\alpha_1>0$ is given by (\ref{ca2}))},\end{array}
\end{equation}
whenever $n$ is large enough so that $2\le e^{2\lambda^*_{L_n}M_{L_n}L_n}$ (remember that $\lambda^*_{L_n}\to+\infty$ as $n\to+\infty$, by assumption). Meanwhile, for all $n\in\N$,
\begin{equation}\label{E(L)}
|C(n)|\le\int_{\theta_n}^{\theta_n+M_{L_n}\,L_n}\left|\mu(\frac{x}{L_n})\right| e^{2\lambda^{*}_{L_n} x}\psi_{n}^{2}(x)dx
\leq\mu_{\infty}\times\frac{\psi_n^2(\theta_n)}{2\lambda^*_{L_n}}\times e^{2\lambda^*_{L_n}(\theta_n+M_{L_n}L_n)},
\end{equation}
where $\mu_{\infty}=\max_{x\in\R}|\mu(x)|$. On the other hand, (\ref{var formula}) and (\ref{c*L is bounded}) yield
$$k_n\le 2\sqrt{a_M\mu_M}\times\lambda^*_{L_n}$$
for all $n\in\N$, whence
\begin{equation}\label{rhs}\baa{rcl}
k_n\displaystyle{\int_{\theta_n}^{\theta_n+M_{L_n}\,L_n}\varphi_{n}^{2}} & = & k_n\displaystyle{\int_{\theta_n}^{\theta_n+M_{L_n}\,L_n}e^{2\lambda^*_{L_n}x}\psi_{n}^{2}}(x)dx\\
& \le & \sqrt{a_M\mu_M}\times\psi_n^2(\theta_n)\times e^{2\lambda^*_{L_n}(\theta_n+M_{L_n}L_n)}.\eaa
\end{equation}
Now, the term $B(n)$ can be estimated as follows
\begin{eqnarray*}\begin{array}{ll}
B(n)&=\displaystyle{\sum_{j=0}^{M_{L_n}-1}\int_{\theta_n+jL_n}^{\theta_n+(j+1)L_n}a_{L_n}e^{2\lambda^{*}_{L_n}x}\left(\psi'_{n}(x)+\lambda^{*}_{L_n}\psi_{n}(x)\right)^{2}dx}\\
&\leq\displaystyle{\sum_{j=0}^{M_{L_n}-1}\alpha_{2}\ e^{2\lambda^{*}_{L_n}(\theta_n+(j+1)L_n)}\int_{\theta_n+jL_n}^{\theta_n+(j+1)L_n}\left(\psi'_{n}(x)+\lambda^{*}_{L_n}\psi_{n}(x)\right)^{2}dx}~\hbox{[from (\ref{ca2})]}\\
&=\displaystyle{\sum_{j=0}^{M_{L_n}-1}\alpha_{2}\ e^{2\lambda^{*}_{L_n}(\theta_n+(j+1)L_n)}\int_{0}^{L_n}\left(\psi'_{n}(x)+\lambda^{*}_{L_n}\psi_{n}(x)\right)^{2}dx}
\end{array}
\end{eqnarray*}
since $\psi_n$ is $L_n$-periodic. One has
$$\begin{array}{ll}
\displaystyle{\int_{0}^{L_n}\left(\psi'_{n}(x)+\lambda^{*}_{L_n}\psi_{n}(x)\right)^{2}dx}&\displaystyle{\leq \;\psi_{n}^{2}(\theta_n)\int_{0}^{L_n}\left(\frac{\psi_{n}'(x)}{\psi_{n}(x)}+\lambda^{*}_{L_n}\right)^{2}dx}.\end{array}$$
We refer now to equation (\ref{elliptic equation L periodic}). Taking $\lambda=\lambda^{*}_{L_n},$ dividing this equation (\ref{elliptic equation L periodic}) by the $L_n$-periodic function $\psi_{n}$ and then integrating by parts over the interval $[0,L_n],$ we get
$$\int_{0}^{L_n}a_{L_n}\left(\frac{\psi'_{n}}{\psi_{n}}\right)^{2}+2\lambda^{*}_{L_n}\int_{0}^{L_n}a_{L_n}\frac{\psi_{n}'}{\psi_{n}}+{\lambda^{*}_{L_n}}^{2}\int_{0}^{L_n}a_{L_n}+\int_{0}^{L_n}\mu_{L_n}=k_nL_n$$
for all $n\in\N$. Thus,
$$\int_{0}^{L_n}a_{L_n}\left(\frac{\psi_{n}'}{\psi_{n}}+\lambda^{*}_{L_n}\right)^{2}+\underbrace{\int_{0}^{L_n}\mu_{L_n}}_{>0}=k_nL_n\leq 2\sqrt{a_M\mu_M}\times\lambda^{*}_{L_n}\,L_n.$$
Owing to (\ref{ca2}), it follows that
$$\forall\ n\in\N,\ \int_{0}^{L_n}\left(\frac{\psi_{n}'}{\psi_{n}}+\lambda^{*}_{L_n}\right)^{2}\leq\frac{2\sqrt{a_M\mu_M}}{\alpha_1}\times\lambda^{*}_{L_n}\,L_n.$$
Putting the above result into $B(n),$ we obtain, for all $n\in\N$,
\begin{equation}\label{B(L)}
\begin{array}{ll}
B(n)&\leq \displaystyle{\frac{2\alpha_2\sqrt{a_M\mu_M}}{\alpha_1}\times\lambda^{*}_{L_n}\,L_n\psi_{n}^{2}(\theta_n)\sum_{j=0}^{M_{L_n}-1}e^{2\lambda^{*}_{L_n}(\theta_n+(j+1)L_n)}}\\\\
&=\displaystyle{\frac{2\alpha_2\sqrt{a_M\mu_M}}{\alpha_1}\times\lambda^{*}_{L_n}\,L_n\psi_{n}^{2}(\theta_n)e^{2\lambda^{*}_{L_n}(\theta_n+L_n)}\times\frac{e^{2\lambda^{*}_{L_n}\,L_nM_{L_n}}-1}{e^{2\lambda^{*}_{L_n}L_n}-1}}\\\\
&\leq\displaystyle{\frac{2\alpha_2\sqrt{a_M\mu_M}}{\alpha_1}\times\psi_{n}^{2}(\theta_n)\times\frac{\lambda^{*}_{L_n}\,L_ne^{2\lambda^{*}_{L_n}L_n}}{e^{2\lambda^{*}_{L_n}L_n}-1}\times e^{2\lambda^{*}_{L_n}(\theta_n+M_{L_n}\,L_n)}}\\\\
&\leq \displaystyle{\beta\times\psi_{n}^{2}(\theta_n)e^{2\lambda^{*}_{L_n}(\theta_n+M_{L_n}\,L_n)}\times\left(\lambda^{*}_{L_n}L_n+1\right)},
\end{array}
\end{equation}
where $\beta=\left(2\alpha_2\sqrt{a_M\mu_M}/\alpha_1\right)\times C$ and $C$ is a positive constant such that
$$\forall x\geq0,~~\frac{xe^{2x}}{e^{2x}-1}\leq C\times(x+1).$$\par
Lastly, let us rewrite equation (\ref{int by parts over theta_L}) as
$$\forall\ n\in\N,~~A(n)+C(n)-k_n\int_{\theta_n}^{\theta_n+M_{L_n}L_n}\varphi_{n}^{2}=B(n).$$
Together with (\ref{A(L)}), (\ref{E(L)}), (\ref{rhs}) and (\ref{B(L)}), one concludes that there exists $n_0\in\N$ such that for $n\ge n_0$,
\begin{equation}\label{mix all}
\begin{array}{l}
\displaystyle{\frac{\alpha_{1}}{2}\times\lambda^{*}_{L_n}\psi^{2}_{n}(\theta_n)e^{2\lambda^{*}_{L_n}(\theta_n+M_{L_n}\,L_n)}}\ -\ \mu_{\infty}\times\displaystyle{\frac{\psi_n^2(\theta_n)}{2\lambda^*_{L_n}}}\times e^{2\lambda^*_{L_n}(\theta_n+M_{L_n}L_n)}\\
-\ \sqrt{a_M\mu_M}\times\psi_n^2(\theta_n)e^{2\lambda^{*}_{L_n}(\theta_n+M_{L_n}\,L_n)}\\
\\
\qquad\qquad\leq\ \displaystyle{\beta\times\psi_{n}^{2}(\theta_n)e^{2\lambda^{*}_{L_n}(\theta_n+M_{L_n}\,L_n)}\times\left(\lambda^{*}_{L_n}L_n+1\right)}.
\end{array}
\end{equation}
Divide (\ref{mix all}) by $\lambda^{*}_{L_n}\psi^{2}_{n}(\theta_n)e^{2\lambda^{*}_{L_n}(\theta_n+M_{L_n}\,L_n)}.$ Then
$$\forall\ n\ge n_0,~~\frac{\alpha_{1}}{2}-\frac{\mu_{\infty}}{2(\lambda^{*}_{L_n})^2}-\frac{\sqrt{a_M\mu_M}}{\lambda^{*}_{L_n}}\leq\beta\times\left(L_n+\frac{1}{\lambda^{*}_{L_n}}\right).$$
Passing to the limit as $n\to+\infty$, one has $L_n\to 0^+$ and $\lambda^*_{L_n}\to+\infty$, whence $\alpha_1\le 0$, which is impossible.\par
Therefore the assumption that $\lambda_{L_n}^{*}\rightarrow+\infty$ as $L_n\rightarrow0^+$ is false and consequently the family $\left(\lambda^{*}_{L}\right)_L$ is bounded from above by some positive $\overline{\lambda}>0$ whenever $L$ is small (i.e. $0<L\leq L_0$). This completes the proof of Lemma \ref{lemma 1}.\hfill{$\Box$}

\begin{remark}\label{lambda*}
{\rm From Theorem \ref{limit as L tends to 0}, one concludes that the map
$(0,+\infty)\ni L\mapsto c^{*}_{L}$ can be extended by continuity to the right at $L=0^+.$ Furthermore, for any sequence $(L_n)_n$ of positive numbers such that
$L_n\rightarrow0^+$ as $n\rightarrow+\infty$, one claims that the positive numbers $\lambda^*_{L_n}$ given in (\ref{var formula}) converge to $\sqrt{<\!a\!>_H^{-1}<\!\mu\!>_A}=\sqrt{<\!a^{-1}\!>_A<\!\mu\!>_A}$ as $n\to+\infty$. Indeed
$$\displaystyle{\forall\ n\in\mathbb{N},~~c^{*}_{L_n}=\frac{k(\lambda^{*}_{L_n},L_n)}{\lambda^{*}_{L_n}}}$$
and Lemma \ref{lemma 1} implies that, up to extraction of a subsequence, $\lambda^{*}_{L_n}\rightarrow\lambda^{*}>0.$ Passing to the limit as $n\rightarrow+\infty$ in the above equation and due (\ref{msh}) together with Step~2 of the proof of Theorem \ref{limit as L tends to 0}, one gets
$$2\sqrt{<\!a\!>_H<\!\mu\!>_A}=\frac{\overline{k}(\lambda^{*})}{\lambda^{*}}=\lambda^{*}<\!a\!>_H+\frac{<\!\mu\!>_A}{\lambda^{*}},$$
whence $\lambda^*=\sqrt{<\!a\!>_H^{-1}<\!\mu\!>_A}$. Since the limit does not depend on any subsequence, one concludes that the limit of $\lambda^{*}_{L},$ as $L\rightarrow0^{+},$ exits and
$$\lim_{L\rightarrow0^+}\lambda^{*}_{L}=\sqrt{<\!a\!>_H^{-1}<\!\mu\!>_A}=\sqrt{<\!a^{-1}\!>_A<\!\mu\!>_A}.$$}
\end{remark}
\vskip 0.2cm

\noindent\textbf{The sharp lower bound of $\liminf_{L\rightarrow0^+}c^{*}_L$ from the homogenized equation}. In the following, we are going to derive the homogenized equation of (\ref{eqevo}), which will lead to the sharp lower bound of $\liminf_{L\rightarrow0^+}c^{*}_L$. However, to furnish this goal we will only consider for the sake of simplicity a particular type of nonlinearities among those satisfying (\ref{cf2}). In fact, the following ideas can be generalized to a wider family of nonlinearities which satisfy (\ref{cf2}), but the proof requires technical extra-arguments which will be the purpose of a forthcoming paper.\par
For each $L>0$, let $u_L$ be a pulsating travelling front with minimal speed $c^*_L$ for the reaction-diffusion equation
\begin{equation}\label{particular}
\left\{\begin{array}{l}
\ds{\frac{\partial u_L}{\partial t}}=\ds{\frac{\partial}{\partial x}\left(a_L(x)\frac{\partial \,u_L}{\partial x}\right)+\mu(\frac{x}{L})g(u_L),\quad t\in\R,\;x\in\R,}\vspace{3 pt}\\
\ds{\forall (t,x)\in\R\times\R,~ 0<u_L(t+\frac{L}{c^*_L},x)=u_L(t,x+L)<1,}\vspace{3 pt}\\
\ds{\lim_{x\rightarrow-\infty}u_L(t,x)=0\hbox{ and }
\lim_{x\rightarrow+\infty}u_L(t,x)=1,}\end{array}
\right.
\end{equation}
where $a_L(x)=a(x/L)$, $a$ is a $C^{2,\delta}(\R)$ $1$-periodic function satisfying (\ref{ca2}), $\mu$ is a $C^{1,\delta}(\R)$ positive 1-periodic function  and $g$ is a $C^2(\R_+)$ function such that $g(0)=g(1)=0$ and $u\mapsto g(u)/u$ is decreasing in $(0,+\infty)$. Up to a shift in time, one can assume that
\begin{equation}\label{normalizing}
\ds{\forall L>0, ~~\int\!\!\!\int_{(0,1)\times(0,1)}u_L(t,x)\ \!dt\ \!dx=\frac{1}{2}.}
\end{equation}
For each $L>0$, set $f_L(x,u):=f(x/L,u)=\mu(x/L)g(u)$. In this setting, there holds $p_L\equiv 1$. From standard parabolic estimates, each function $u_L$ is (at least) of class $C^2(\R\times\R)$. Denote
$$v_L(t,x)=a_L(x)\frac{\partial u_L}{\partial x}(t,x)\;\hbox{ and }\;w_L(t,x)=\frac{\partial u_L}{\partial t}(t,x)\; \hbox{ in }\;\R\times\R.$$
As already underlined, it follows from \cite{bh} that $w_L=\frac{\partial u_L}{\partial t}>0$ in $\R\times\R$ for each $L>0$. Under the notations of the beginning of this section, it follows from (\ref{ca2}) and (\ref{elliptic equation L periodic}) that $k(\lambda,L)\ge\lambda^2\alpha_1+\mu_m$ for all $L>0$ and $\lambda\in\R$, where $\mu_m=\min_{\R}\mu>0$. Hence, $c^*_L\ge2\sqrt{\alpha_1\mu_m}$ for each $L>0$ and $\liminf_{L\to 0^+}c^*_L\ge2\sqrt{\alpha_1\mu_m}>0$.\par
We shall now establish some estimates for the functions $u_L$, $v_L$ and $w_L$ which are independent of $L$, in order to pass to the limit as $L\to 0^+$. Notice first that standard parabolic estimates and the $(t,x)$-periodicity satisfied by the functions $u_L$ imply that, for each $L>0,$ $u_L(-\infty,x)=0$ and $u_L(+\infty,x)=1$ in $C^2_{loc}(\R)$, and $w_L(\pm\infty,x)=0$ in $C^1_{loc}(\R)$.\par
Let $k\in\N\backslash\{0\}$ be given. Integrating the first equation of (\ref{particular}) by parts over $\R\times (-kL,kL),$ one obtains
\be\label{intf}
\ds{\forall L>0,~\int\!\!\!\int_{\R\times(-kL,kL)}f(\frac{x}{L},u_L)\ \!dt\ \!dx=2kL.}
\ee
Multiplying the first equation of (\ref{particular}) by $u_L$ and integrating by parts over $\R\times (-kL,kL),$ one then gets
\begin{equation}\label{inte}
\forall L>0,~kL=-\int\!\!\!\int_{\R\times(-kL,kL)}a_L(x)\left(\frac{\partial u_L}{\partial x}\right)^2\ \!dt\ \!dx+\int\!\!\!\int_{\R\times(-kL,kL)}f(\frac{x}{L},u_L)u_L\ \!dt\ \!dx.
\end{equation}
Notice that the last integral in (\ref{inte}) converges because of (\ref{intf}) and $0\leq f (x/L,u_L)u_L\leq f(x/L,u_L)$. Together with (\ref{ca2}), one concludes that for each $L>0,$ the first integral in~(\ref{inte}) converges and
$$\forall L>0,~\int\!\!\!\int_{\R\times(-kL,kL)}\left(\frac{\partial u_L}{\partial x}\right)^2\ \!dt\ \!dx\leq \frac{kL}{\alpha_1}.$$
Multiply the first equation of (\ref{particular}) by $\frac{\partial u_L}{\partial t}$ and integrate by parts over $\R\times(-kL,kL).$ Since
$$\ds{\int\!\!\!\int_{\R\times(-kL,kL)}\frac{\partial}{\partial x}\left(a_L(x)\frac{\partial u_L}{\partial x}\right)\frac{\partial u_L}{\partial t}}=\ds{-\frac{1}{2}\int_{\R\times(-kL,kL)} \frac{\partial}{\partial t}\left(a_L(x)\left(\frac{\partial u_L}{\partial x}\right)^2\right)}=0,$$
one obtains that
\be\label{int2}
\forall L>0,~\int\!\!\!\int_{\R\times(-kL,kL)}\left(\frac{\partial u_L}{\partial t}\right)^2\ \!dt\ \!dx=\int_{(-kL,kL)}F(\frac{x}{L},1)dx=2kL\times\int_0^1\mu\times\int_0^1g,
\ee
where $F(y,s)=\int_{0}^{s}f(y,\tau)d\tau$. It follows from the above estimates that for each compact subset $K$ of $\R,$
\begin{equation}\label{est H1 for u_L}
\forall\ \! 0<L<1,\;\ds{\int\!\!\!\int_{\R\times K}}\left [\left(\frac{\partial u_L}{\partial t}\right)^2+\left(\frac{\partial u_L}{\partial x}\right)^2\right]dt\ \!dx\leq C(K),
\end{equation}
where $C(K)$ is a positive constant depending only on $K.$\par
In particular, for each compact $K$ of  $\R$ and for each $L>0$, $||w_L||_{L^2(\R\times K)}\leq \sqrt{C(K)}$. Now, differentiate the first equation of (\ref{particular}) with respect to $t$ (actually, from the regularity of $f$, the function $w_L$ is of class $C^2$ with respect to $x$). There holds
$$\ds{\frac{\partial w_L}{\partial t}}=\ds{\frac{\partial}{\partial x}\left(a_L(x)\frac{\partial \,w_L}{\partial x}\right)+\mu(\frac{x}{L})g'(u_L)w_L\ \hbox{ in }\R\times\R.}$$
Multiply the above equation by $w_L$ and integrate by parts over $\R\times(-kL,kL)$. From (\ref{ca2}) and (\ref{int2}), it follows that
$$\ds{\int\!\!\!\int_{\R\times(-kL,kL)}\left(\frac{\partial w_L}{\partial x}\right)^2dtdx\leq \frac{2kL\eta}{\alpha_1}}$$
where $\eta$ is the positive constant defined by $$\ds{\eta=\max_{x\in\R}\mu(x)\max_{u\in[0,1]}|g'(u)|\,\max_{x\in\R}|F(x,1)|\geq\frac{1}{2kL}\int\!\!\!\int_{\R\times(-kL,kL)}\mu(\frac{x}{L})g'(u_L)w_L^2\ \!dt\ \!dx}>0.$$
Then, for each compact $K\subset\R,$ there exists a constant $C'(K)>0$ depending only on $K$ such that
\begin{equation}\label{est H1 for w_L}
\forall\ \! 0<L<1,\;\ds{\int\!\!\!\int_{\R\times K}}\left(\frac{\partial w_L}{\partial x}\right)^2dt\ \!dx\leq C'(K).
\end{equation}\par
Let $(L_n)_{n\in\N}$ be a sequence of real numbers in $(0,1)$ such that $L_n\to 0$ and $c^*_{L_n}\to\liminf_{L\to0^+}c^*_L>0$ as $n\to+\infty$. It follows from (\ref{est H1 for u_L}) and the bounds $0<u_{L_n}<1$ that there exists $u_0$ in $H^{1}_{loc}(\R\times \R)$ such that, up to extraction of a subsequence, $u_{L_n}\rightarrow u_0$ strongly in $L^2_{loc}(\R\times\R)$ and almost everywhere
in $\R\times \R$, and
$$\left(\frac{\partial u_{L_n}}{\partial t},\frac{\partial u_{L_n}}{\partial x}\right)\rightharpoonup\left(\frac{\partial u_0}{\partial t}, \frac{\partial u_0}{\partial x}\right)\hbox{ weakly in }L^2_{loc}(\R\times \R)\hbox{ as }n\rightarrow+\infty.$$
Remember that $v_{L_n}=a_{L_n}\frac{\partial u_{L_n}}{\partial x}$  and $0<\alpha_1\leq a_{L_n}\leq\alpha_2 $ for each $n\in\N$.  Thus, (\ref{est H1 for u_L}) yields that for each compact $K$ of $\R$ and for each $n\in\N$, $\ds{||v_{L_n}||_{L^2(\R\times K)}\leq\alpha_2C(K).}$ Furthermore,~(\ref{particular}) implies that
$$\forall n\in\N,~\ds{\frac{\partial v_{L_n}}{\partial x}=\frac{\partial u_{L_n}}{\partial t}-f(\frac{x}{L_n},u_{L_n})}\hbox{ in }\R\times\R,$$
while $0\le f(x/L_n,u_{L_n}(t,x))\leq \kappa$ in $\R\times\R$ where $\kappa=\max_{\R}\mu\times\max_{[0,1]}g>0$ is independent of~$n.$  Together with (\ref{est H1 for u_L}), one concludes that the sequence
$(\frac{\partial v_{L_n}}{\partial x})_{n\in\N}$ is bounded in $L^2_{loc}(\R\times\R).$ On the other hand, $\frac{\partial v_{L_n}}{\partial t}=a_{L_n}\frac{\partial^2u_{L_n}}{\partial t\partial x}=a_{L_n}\frac{\partial w_L }{\partial x}$. Owing to (\ref{ca2}) and (\ref{est H1 for w_L}), the sequence $(\frac{\partial v_{L_n}}{\partial t})_{n\in\N}$ is bounded in
$L^2_{loc}(\R\times\R).$ Consequently, up to extraction of another subsequence, there exists $v_0\in H^{1}_{loc}(\R\times\R)$ such that $v_{L_n}\to v_0$ strongly in $L^2_{loc}(\R\times\R)$ and
$$\left(\frac{\partial v_{L_n}}{\partial t},\frac{\partial v_{L_n}}{\partial x}\right)\rightharpoonup\left(\frac{\partial v_{0}}{\partial t},\frac{\partial v_{0}}{\partial x}\right)\hbox{ weakly in }L^2_{loc}(\R\times\R)\hbox{ as }n\to+\infty.$$
However, $a_{L_n}^{-1}\rightharpoonup\,<a^{-1}>_A=<a>_H^{-1}\hbox{ in }L^{\infty}(\R)$ weak-$*$ as $n\to+\infty$. Thus,
$$\frac{\partial u_{L_n}}{\partial x}=\frac{v_{L_n}}{a_{L_n}}\rightharpoonup\frac{v_0}{<a>_H}\hbox{ weakly in }L^2_{loc}(\R\times\R)\hbox{ as }n\to+\infty.$$
By uniqueness of the limit, one gets  $v_0=<a>_H\frac{\partial u_0}{\partial x}$. Passing to the limit as $n\to+\infty$ in the first equation of (\ref{particular}) with $L=L_n$ implies that $u_0$ is a weak solution of the equation
 $$\frac{\partial u_0}{\partial t}=\frac{\partial v_0}{\partial x}+<\mu>_Ag(u_0)=<a>_H\frac{\partial^2 u_0}{\partial x^2}+<\mu>_Ag(u_0)\hbox{ in }\mathcal{D}'(\R\times\R).$$
From parabolic regularity, the function $u_0$ is then a classical solution of the homogenous equation
$$\frac{\partial u_0}{\partial t}=<a>_H \frac{\partial^{2}u_0}{\partial x ^2}+<\mu>_A \,g(u_0)\hbox{ in }\R\times\R,$$
such that $0\le u_0\le 1$ and $\frac{\partial u_0}{\partial t}\geq0$ in $\R\times\R$. Lastly, $\int\!\!\!\int_{(0,1)^2}u_0(t,x)\ \!dt\ \!dx=\frac{1}{2}$ from (\ref{normalizing}). On the other hand, it follows from the second equation of (\ref{particular}) and (\ref{est H1 for u_L}) that
$$\forall \gamma\in \R,~~u_0(t+\frac{\gamma}{c},x)=u_0(t,x+\gamma)\ \hbox{ in }\R\times\R,$$
where $c=\liminf_{L\to0^+}c^*_L=\lim_{n\to+\infty}c^*_{L_n}>0$. In other words, $u_0(t,x)=U_0(x+ct)$, where $U_0$ is a classical solution of the equation
\begin{equation}\label{hom U_0}
\ds{cU'_0=<a>_HU''_0+<\mu>_Ag(U_0),\ \ 0\le U_0\le 1~\hbox{ in }\R}
\end{equation}
that satisfies $U_0'\ge 0$ in $\R$ and
$$\int_0^1\left(\int_{cs}^{cs+1}U_0\right)ds=\frac{1}{2}.$$
Standard elliptic estimates imply that $U_0$ converges as $s\to\pm\infty$ in $C^2_{loc}(\R)$ to two constants $U^{\pm}_0\in[0,1]$ such that $<\mu>_Ag(U^{\pm}_0)=0$, that is $g(U^{\pm}_0)=0$. The monotonicity of $U_0$ and the assumption on $g$ imply that $U^-_0=0$ and $U^+_0=1.$ In other words, $U_0$ is a usual travelling front for the homogenized equation (\ref{hom U_0}) with speed $c$ and limiting conditions $0$ and $1$ at infinity. Since the minimal speed for this problem is equal to $\ds{2\sqrt{<a>_H<\mu>_A}}$, one concludes that $$\liminf_{L\rightarrow0^+}c^*_L=c\geq \ds{2\sqrt{<a>_H<\mu>_A}}.$$

%%%%%%%%%%%%%%%%%%%%%%%%%%%%%%%%%%%%%%%%
%%%%%%%%%%%%%%%%%%%%%%%%%%%%%%%%%%%%%%%%

\section{Monotonicity of the minimal speeds $c^*_L$ near the homogenization limit}\label{variation}

This section is devoted to the proof of Theorem \ref{c*L is increasing}. Before going further in the proof, we recall that for each $L>0,$ the minimal speed $c^{*}_{L}$ is given by the variational formula $$c^{*}_L=\min_{\lambda>0}\frac{k(\lambda,L)}{\lambda}=\frac{k(\lambda^{*}_L,L)}{\lambda^{*}_L},$$
where $\lambda^*_L>0$ and $k(\lambda,L)$ is the principal eigenvalue of the elliptic equation (\ref{elliptic equation L periodic}). Notice that $k(\lambda,L)$ can be defined for all $\lambda\in\R$ and $L>0$.\hfill\break\par

{\underline{Step 1: properties of $k(\lambda,L)$ and definition of $\tilde{k}(\lambda,L)$.} The principal eigenfunction $\psi_{\lambda,L}$ of (\ref{elliptic equation L periodic}) is $L$-periodic, positive and unique up to multiplication. Denote
$$\phi_{\lambda,L}(x)=\psi_{\lambda,L}(Lx)$$
for all $L>0$, $\lambda\in\R$ and $x\in\R$. Each function $\phi_{\lambda,L}$ is $1$-periodic, positive and it is the principal eigenfunction of
$$(a\phi_{\lambda,L}')'+2L\lambda a\phi_{\lambda,L}'+L\lambda a'\phi_{\lambda,L}+L^2\lambda^2a\phi_{\lambda,L}+L^2\mu\phi_{\lambda,L}=L^2k(\lambda,L)\phi_{\lambda,L},$$
associated to the principal eigenvalue $L^2k(\lambda,L)$. But the above problem can be defined for all $\lambda\in\R$ and $L\in\R$. That is, for each $(\lambda,L)\in\R^2$, there exists a unique principal eigenvalue $\tilde{k}(\lambda,L)$ and a unique (up to multiplication) principal eigenfunction $\tilde{\phi}(\lambda,L)$ of
\be\label{tildek}
(a\tilde{\phi}_{\lambda,L}')'+2L\lambda a\tilde{\phi}_{\lambda,L}'+L\lambda a'\tilde{\phi}_{\lambda,L}+L^2\lambda^2a\tilde{\phi}_{\lambda,L}+L^2\mu\tilde{\phi}_{\lambda,L}=\tilde{k}(\lambda,L)\tilde{\phi}_{\lambda,L}.
\ee
Furthermore, $\tilde{\phi}_{\lambda,L}$ is $1$-periodic, positive and it can be normalized so that
\be\label{normalization}
\int_0^1\tilde{\phi}_{\lambda,L}^2(x)dx=1
\ee
for all $(\lambda,L)\in\R^2$. By uniqueness of the principal eigenelements, it follows that
$$\forall\ L>0,\ \forall\ \lambda\in\R,\quad\tilde{k}(\lambda,L)=L^2k(\lambda,L)$$
and $\tilde{\phi}_{\lambda,L}$ and $\phi_{\lambda,L}$ are equal up to multiplication by positive constants for each $L>0$ and $\lambda\in\R$.\par
Some useful properties of $k(\lambda,L)$ as $L\to 0^+$ shall now be derived from the study the function $\tilde{k}$. Notice first that, since the coefficients of the left-hand side of (\ref{tildek}) are analytic in $(\lambda,L)$, the function $\tilde{k}$ is analytic, and from the normalization (\ref{normalization}), the functions $\tilde{\phi}_{\lambda,L}$ also depend analytically in $H^2_{loc}(\R)$ on the parameters $\lambda$ and $L$ (see \cite{cv,kato}). In particular, the function $k$ is analytic in $\R\times(0,+\infty)$. Observe also that
$$\tilde{k}(\lambda,0)=0\hbox{ and }\tilde{\phi}_{\lambda,0}=1\hbox{ for all }\lambda\in\R.$$
Lastly, when $\lambda$ is changed into $-\lambda$ or when $L$ is changed into $-L$, then the operator in (\ref{tildek}) is changed into its adjoint. But since the principal eigenvalues of the operator and its adjoint are identical, it follows that
$$\forall\ (\lambda,L)\in\R^2,\quad \tilde{k}(\lambda,L)=\tilde{k}(\lambda,-L)=\tilde{k}(-\lambda,L).$$
In particular, it follows that
\be\label{partialtildek}\forall\ (i,j)\in\N^2,\quad
\frac{\partial^i\tilde{k}}{\partial\lambda^i}(\lambda,0)=\frac{\partial^i\partial^{2j+1}\tilde{k}}{\partial\lambda^i\partial L^{2j+1}}(\lambda,0)=0.
\ee\par
Therefore, for all $\overline{\lambda}\in\R$,
$$k(\lambda,L)=\frac{\tilde{k}(\lambda,L)}{L^2}\to\frac{1}{2}\times\frac{\partial^2\tilde{k}}{\partial L^2}(\overline{\lambda},0)\hbox{ as }(\lambda,L)\to(\overline{\lambda},0^+).$$
But since this limit is equal to $\overline{k}(\overline{\lambda})=\overline{\lambda}^2<\!a\!>_H+<\!\mu\!>_A$ from Step~2 of the proof of Theorem~\ref{thm: c*L--->** as L--->0}, one then gets that
\be\label{tildek02}
\frac{1}{2}\times\frac{\partial^2\tilde{k}}{\partial L^2}(\overline{\lambda},0)=\overline{\lambda}^2<\!a\!>_H+<\!\mu\!>_A\hbox{ for all }\overline{\lambda}\in\R.
\ee\par
It also follows from (\ref{partialtildek}) that
\be\label{tildek22}
\frac{\partial^2k}{\partial\lambda^2}(\lambda,L)=\frac{1}{L^2}\times\frac{\partial^2\tilde{k}}{\partial\lambda^2}(\lambda,L)\to\frac{1}{2}\times\frac{\partial^4\tilde{k}}{\partial\lambda^2\partial L^2}(\overline{\lambda},0)\hbox{ as }(\lambda,L)\to(\overline{\lambda},0^+).
\ee
From (\ref{tildek02}) and (\ref{tildek22}), one deduces that
\be\label{k20}\frac{\partial^2k}{\partial\lambda^2}(\lambda,L)\to2<\!a\!>_H\ \ >0\hbox{ as }(\lambda,L)\to(\overline{\lambda},0^+).
\ee\par
Similarly, as $(\lambda,L)\to(\overline{\lambda},0^+)$,
\be\label{partialk}\left\{\baa{rcl}
\displaystyle{\frac{\partial k}{\partial L}}(\lambda,L)=\displaystyle{\frac{\partial}{\partial L}}\left(\displaystyle{\frac{\tilde{k}(\lambda,L)}{L^2}}\right) & \to & \displaystyle{\frac{1}{6}}\times\displaystyle{\frac{\partial^3\tilde{k}}{\partial L^3}}(\overline{\lambda},0)=0\\
\\
\displaystyle{\frac{\partial^2k}{\partial\lambda\partial L}}(\lambda,L)=\displaystyle{\frac{\partial}{\partial L}}\left(\displaystyle{\frac{1}{L^2}}\times\displaystyle{\frac{\partial\tilde{k}}{\partial\lambda}}(\lambda,L)\right) & \to & \displaystyle{\frac{1}{6}}\times\displaystyle{\frac{\partial^4\tilde{k}}{\partial\lambda\partial L^3}}(\overline{\lambda},0)=0\\
\\
\displaystyle{\frac{\partial^2k}{\partial L^2}}(\lambda,L)=\displaystyle{\frac{\partial^2}{\partial L^2}}\left(\displaystyle{\frac{\tilde{k}(\lambda,L)}{L^2}}\right) & \to & \displaystyle{\frac{1}{12}}\times\displaystyle{\frac{\partial^4\tilde{k}}{\partial L^4}}(\overline{\lambda},0)\eaa\right.
\ee

\begin{remark}{\rm As a byproduct of the fact that $\tilde{k}$ and $k$ are even in $\lambda$, it follows that the minimal speed of pulsating fronts propagating from right to left (as in Definition~\ref{pulsating traveling fronts}) is the same as that of fronts propagating from left to right.}
\end{remark}

{\underline{Step 2: properties of $c^*_L$ and $\lambda^*_L$ in the neighbourhood of $L=0^+$.}} Let us first prove that, for each fixed $L>0$, the positive real number $\lambda^*_L>0$ given in (\ref{var formula}) is unique. Indeed, if there are $0<\lambda_1<\lambda_2$ such that
$$c^*_L=\frac{k(\lambda_1,L)}{\lambda_1}=\frac{k(\lambda_2,L)}{\lambda_2}=\min_{\lambda>0}\frac{k(\lambda,L)}{\lambda},$$
then $k(\lambda,L)=c^*_L\lambda$ for all $\lambda\in[\lambda_1,\lambda_2]$ since $k$ is convex with respect to $\lambda$. Then $k(\lambda,L)=c^*_L\lambda$ for all $\lambda\in\R$ by analyticity of the map $\R\ni\lambda\mapsto k(\lambda,L)$. But $k(0,L)=-\rho_{1,L}>0$, which gives a contradiction. Therefore, for each $L>0$, $\lambda^*_L$ is the unique minimum of the map $(0,+\infty)\ni\lambda\mapsto k(\lambda,L)/\lambda$.\par
Furthermore, we claim that $L\mapsto\lambda^*_L$ and $L\mapsto c^*_L$ are of class $C^{\infty}$ in a right neighbourhood of $L=0$. Indeed, by definition, $\lambda^*_L$ satisfies
\be\label{F}
F(\lambda^*_L,L):=\frac{\partial k}{\partial\lambda}(\lambda^*_L,L)\times\lambda^*_L-k(\lambda^*_L,L)=0.
\ee
The function $(\lambda,L)\mapsto F(\lambda,L)$ is of class $C^{\infty}$ on $\R\times(0,+\infty)$ and $\frac{\partial F}{\partial\lambda}(\lambda,L)=\frac{\partial^2k}{\partial\lambda^2}(\lambda,L)\times\lambda$. But
$$\lambda^*_L\to\lambda^*=\sqrt{<\!a\!>_H^{-1}<\!\mu\!>_A}\ >\ 0\hbox{ as }L\to 0^+$$
from Remark~\ref{lambda*}, and
$$\frac{\partial^2k}{\partial\lambda^2}(\lambda^*_L,L)\ \to\ 2\ <\!a\!>_H\ \ >0\hbox{ as }L\to 0^+$$
from (\ref{k20}). Therefore, from the implicit function theorem, the map $L\mapsto\lambda^*_L$ is of class $C^{\infty}$ in an interval $(0,L_0)$ for some $L_0>0$. As a consequence of formula (\ref{var formula}), the map $L\mapsto c^*_L$ is also of class $C^{\infty}$ on $(0,L_0)$.\par
For each $L\in(0,L_0)$, one has
$$\frac{dc^*_L}{dL}=\left(\frac{1}{\lambda^*_L}\times\frac{\partial k}{\partial\lambda}(\lambda^*_L,L)-\frac{k(\lambda^*_L,L)}{(\lambda^*_L)^2}\right)\times\frac{d\lambda^*_L}{dL}+\frac{1}{\lambda^*_L}\times\frac{\partial k}{\partial L}(\lambda^*_L,L)=\frac{1}{\lambda^*_L}\times\frac{\partial k}{\partial L}(\lambda^*_L,L)$$
by definition of $\lambda^*_L$ and formula (\ref{var formula}). But $\lambda^*_L\to\lambda^*>0$ and $\frac{\partial k}{\partial L}(\lambda^*_L,L)\to0$ as $L\to 0^+$ from (\ref{partialk}). Thus,
$$\frac{dc^*_L}{dL}\to 0\hbox{ as }L\to 0^+.$$
On the other hand, it follows from (\ref{k20}), (\ref{partialk}) and (\ref{F}) that
$$\frac{d\lambda^*_L}{dL}=\frac{1}{\lambda^*_L\times\displaystyle{\frac{\partial^2k}{\partial\lambda^2}}(\lambda^*_L,L)}\times\left(\frac{\partial k}{\partial L}(\lambda^*_L,L)-\lambda^*_L\times\frac{\partial^2k}{\partial\lambda\partial L}(\lambda^*_L,L)\right)\to 0\hbox{ as }L\to 0^+.$$
Therefore,
\be\label{d2}\baa{rcl}
\displaystyle{\frac{d^2c^*_L}{dL^2}} & = & \displaystyle{\frac{d\lambda^*_L}{dL}}\times\left(-\displaystyle{\frac{1}{(\lambda^*_L)^2}}\times\displaystyle{\frac{\partial k}{\partial L}}(\lambda^*_L,L)+\displaystyle{\frac{1}{\lambda^*_L}}\times\displaystyle{\frac{\partial^2k}{\partial\lambda\partial L}}(\lambda^*_L,L)\right)+\displaystyle{\frac{1}{\lambda^*_L}}\times\displaystyle{\frac{\partial^2k}{\partial L^2}}(\lambda^*_L,L)\\
\\
& \to & \displaystyle{\frac{1}{12\lambda^*}}\times\displaystyle{\frac{\partial^4\tilde{k}}{\partial L^4}}(\lambda^*,0)\ \hbox{ as }L\to 0^+,\eaa
\ee
from (\ref{partialk}).\hfill\break\par

{\underline{Step 3: calculation of $\frac{\partial^4\tilde{k}}{\partial L^4}(\lambda^*,0)$.} In this step, we fix $\lambda^*=\sqrt{<\!a\!>_H^{-1}<\!\mu\!>_A}$. Since the functions $\tilde{\phi}_{\lambda^*,L}$ depend analytically on $L\in\R$ in $H^2_{loc}(\R)$, the expansion
$$\tilde{\phi}_{\lambda^*,L}=1+L\phi_1+L^2\phi_2+L^3\phi_3+L^4\phi_4+\ldots$$
is valid in $H^2_{loc}(\R)$ in a neighbourhood of $L=0$, where $1=\tilde{\phi}_{\lambda^*,0}$ and
$$\phi_i=\left.\frac{1}{i\ !}\times\frac{\partial^i\tilde{\phi}_{\lambda^*,L}}{\partial L^i}\right|_{L=0}$$
for each $i\ge 1$. We now put this expansion into
$$(a\tilde{\phi}_{\lambda^*,L}')'+2L\lambda^*a\tilde{\phi}_{\lambda^*,L}'+L\lambda^*a'\tilde{\phi}_{\lambda^*,L}+L^2(\lambda^*)^2a\tilde{\phi}_{\lambda^*,L}+L^2\mu\tilde{\phi}_{\lambda^*,L}=\tilde{k}(\lambda^*,L)\tilde{\phi}_{\lambda^*,L}$$
and remember that
$$\tilde{k}(\lambda^*,0)=\frac{\partial\tilde{k}}{\partial L}(\lambda^*,0)=\frac{\partial^3\tilde{k}}{\partial L^3}(\lambda^*,0)=0$$
and
$$\frac{\partial^2\tilde{k}}{\partial L^2}(\lambda^*,0)=2\times\left[(\lambda^*)^2<\!a\!>_H+<\!\mu\!>_A\right]=4<\!\mu\!>_A$$
from (\ref{partialtildek}) and (\ref{tildek02}). Since both $\tilde{\phi}_{\lambda^*,L}$ and $\tilde{k}(\lambda^*,L)$ depend analytically on $L$, it follows in particular that
\be\label{phii}\left\{\baa{l}
(a\phi_1')'+\lambda^*a'=0,\\
(a\phi_2')'+2\lambda^*a\phi_1'+\lambda^*a'\phi_1+(\lambda^*)^2a+\mu=2<\!\mu\!>_A,\\
(a\phi_3')'+2\lambda^*a\phi_2'+\lambda^*a'\phi_2+(\lambda^*)^2a\phi_1+\mu\phi_1=2<\!\mu\!>_A\phi_1,\\
(a\phi_4')'+2\lambda^*a\phi_3'+\lambda^*a'\phi_3+(\lambda^*)^2a\phi_2+\mu\phi_2=2<\!\mu\!>_A\phi_2+\displaystyle{\frac{1}{24}}\times\displaystyle{\frac{\partial^4\tilde{k}}{\partial L^4}}(\lambda^*,0)\eaa\right.
\ee
in $\R$. Furthermore, each function $\phi_i$ is $1$-periodic and, by differentiating the normalization condition $\|\tilde{\phi}_{\lambda^*,L}\|_{L^2(0,1)}^2=1$ with respect to $L$ at $L=0$, it follows especially that
$$\int_0^1\phi_1=0\hbox{ and }\int_0^1\phi_2=-\frac{1}{2}\int_0^1\phi_1^2.$$\par
It is then found that, for all $x\in\R$,
$$\phi_1(x)=\lambda^*\times\left(-x+<\!a\!>_H\int_0^x\frac{1}{a(y)}dy-\frac{1}{2}-<\!a\!>_H\int_0^1\frac{y}{a(y)}dy\right)$$
and
$$\baa{rcl}
\phi'_2(x) & = & <\!\mu\!>_A\times\left[\displaystyle{\frac{x}{a(x)}}-\displaystyle{\int_0^1}\displaystyle{\frac{y}{a(y)}}dy-\displaystyle{\int_0^x}\displaystyle{\frac{1}{a(y)}}dy-\displaystyle{\frac{<\!a\!>_H}{a(x)}}\displaystyle{\int_0^1}\displaystyle{\frac{y}{a(y)}}dy\right]\\
\\
& & +\displaystyle{\frac{1}{a(x)}}\times\left[<\!a\!>_H\displaystyle{\int_0^1}\left(\displaystyle{\frac{1}{a(y)}}\displaystyle{\int_0^y}\mu(z)dz\right)dy-\displaystyle{\int_0^x}\mu(y)dy\right]+(\lambda^*)^2\left(x+\displaystyle{\frac{1}{2}}\right).\eaa$$
Moreover, it follows from the third equation of (\ref{phii}) that, for all $x\in\R$,
$$\baa{rcl}
a(x)\phi_3'(x) & = & -2\lambda^*\displaystyle{\int_0^x}a(y)\phi'_2(y)dy-\lambda^*\displaystyle{\int_0^x}a'(y)\phi_2(y)dy\\
\\
& & -(\lambda^*)^2\displaystyle{\int_0^x}a(y)\phi_1(y)dy-\displaystyle{\int_0^x}\mu(y)\phi_1(y)dy+2<\!\mu\!>_A\displaystyle{\int_0^x}\phi_1(y)dy+<\!a\!>_Hc,\eaa$$
where
$$\baa{l}
c=\displaystyle{\int_0^1}\left[\displaystyle{\frac{1}{a(y)}}\times\left(2\lambda^*\displaystyle{\int_0^y}a(z)\phi'_2(z)dz+\lambda^*\displaystyle{\int_0^y}a'(z)\phi_2(z)dz+(\lambda^*)^2\displaystyle{\int_0^y}a(z)\phi_1(z)dz\right.\right.\\
\qquad\qquad\qquad\qquad\left.\left.+\displaystyle{\int_0^y}\mu(z)\phi_1(z)dz-2<\!\mu\!>_A\displaystyle{\int_0^y}\phi_1(z)dz\right)\right]dy\eaa$$\par
On the other hand, by integrating the fourth equation of (\ref{phii}) over the interval $[0,1]$, one gets that
\be\label{fourth}
\displaystyle{\frac{1}{24}}\times\displaystyle{\frac{\partial^4\tilde{k}}{\partial L^4}}(\lambda^*,0)=-2<\!\mu\!>_A\int_0^1\phi_2+\lambda^*\int_0^1a\phi_3'+(\lambda^*)^2\int_0^1a\phi_2+\int_0^1\mu\phi_2.
\ee
Now, put all the previous calculations into (\ref{fourth}). After a lengthy sequence of integrations by parts, it is finally found that
$$\displaystyle{\frac{1}{24}}\times\displaystyle{\frac{\partial^4\tilde{k}}{\partial L^4}}(\lambda^*,0)=\int_0^1\frac{A(x)^2}{a(x)}dx\ -\ <\!a\!>_H\left(\int_0^1\frac{A(x)}{a(x)}dx\right)^2,$$
where
$$A(x)=\int_0^x\mu(y)dy\ +\ <\!\mu\!>_A<\!a\!>_H\int_0^x\frac{1}{a(y)}dy\ -\ 2<\!\mu\!>_A\ x.$$\par
From (\ref{d2}), it follows that
$$\displaystyle{\frac{d^2c^*_L}{dL^2}}\to\gamma:=2\sqrt{<\!a\!>_H<\!\mu\!>_A^{-1}}\times\left[\int_0^1\frac{A(x)^2}{a(x)}dx\ -\ <\!a\!>_H\left(\int_0^1\frac{A(x)}{a(x)}dx\right)^2\right]\hbox{ as }L\to 0^+.$$
Cauchy-Schwarz inequality yields $\gamma\ge 0$. Furthermore, $\gamma=0$ if and only if $A$ is constant. But since $A(0)=0$, the condition $\gamma=0$ is equivalent to $A'(x)=0$ for all $x$, which means that
$$\frac{\mu(x)}{<\!\mu\!>_A}+\frac{<\!a\!>_H}{a(x)}=2\ \hbox{ for all }x\in\R.$$
In particular, if $\mu$ is constant and $a$ is not constant (resp. if  $a$ is constant and $\mu$ is not constant), then this condition is not satisfied, whence $\lim_{L\to 0^+}\frac{d^2c^*_L}{dL^2}>0$ in this case. That completes the proofs of Theorem~\ref{c*L is increasing} and Corollary~\ref{corollaire}.\hfill{$\Box$}

\begin{remark}{\rm In the case when $<\!\mu\!>_A=0$ and $\mu\not\equiv 0$, then $\rho_{1,L}<0$ for each $L>0$, and the minimal speed $c^*_L$ of pulsating traveling fronts is well-defined and it is still positive. From the arguments developed in this section and in the previous one, one can check that, in this case,
$$c^*_L\to 0^+,\quad\lambda^*_L\to 0^+,\quad\frac{d\lambda^*_L}{dL}\to\sqrt{\frac{\beta}{<\!a\!>_H}}>0\ \hbox{ and }\ \frac{dc^*_L}{dL}\to2\sqrt{\beta<\!a\!>_H}>0\hbox{ as }L\to 0^+,$$
where
$$\beta=\int_0^1\frac{A(x)^2}{a(x)}dx\ -\ <\!a\!>_H\left(\int_0^1\frac{A(x)}{a(x)}dx\right)^2>0$$
and $A(x)=\displaystyle{\int_0^x}\mu(y)dy$. Therefore, the speeds $c^*_L$ are increasing in a right neighbourhood of $L=0$ but, in this case, the variation is of the first order. Notice that the formula $\lim_{L\to 0^+}\frac{dc^*_L}{dL}=2\sqrt{\beta<\!a\!>_H}$ is coherent with the numerical calculations done by Kinezaki, Kawasaki and Shigesada in \cite{kks} (see Figure~3b with $<\!\mu\!>_A=0$, that is $A=0$ under the notations of \cite{kks}).}
\end{remark}

%%%%%%%%%%%%%%%%%%%%%%%%%%%%%%%%%%%%%%%%
%%%%%%%%%%%%%%%%%%%%%%%%%%%%%%%%%%%%%%%%

\section{Proof of Theorem \ref{th_fragmentation} \label{proof th frag}}

As in the proofs of the previous theorems, we use the following
formula for the minimal speed: \be c^*_z
=\min_{\lambda>0}\frac{k_z(\lambda)}{\lambda}=\frac{k_z(\lambda^{*}_z)}{\lambda^{*}_z},\label{forc*_frag}\ee
where $k_z(\lambda)$ is defined as the unique real number such
that there exists a positive $L_0$-periodic function $\psi$
satisfying:
\begin{equation}\label{Lcmu}
\psi''+2\lambda\ \psi'
+\lambda^{2}\psi+\mu_z(x)\psi=k_z(\lambda)\psi \hbox{ in }(0,L_0).
\end{equation}

Setting $\varphi(x)=e^{\lambda x}\psi(x)$, the above equation and
periodicity conditions become equivalent to: \be \left\{ \baa{l}
\varphi''+\mu_z(x)\varphi=k_z(\lambda) \varphi \hbox{ in }(0,L_0),\\
\varphi(L_0)=e^{\lambda L_0}\varphi(0), \\
\varphi'(L_0)=e^{\lambda L_0}\varphi'(0),  \eaa \right.
\label{sy1} \ee which therefore admits, for every positive
$\lambda$, a unique solution $(\varphi,k_z(\lambda))$ with
$\varphi>0$ satisfying the normalisation condition $\varphi(0)=1$.

Let $\lambda>0$ be fixed. System (\ref{sy1}), together with the
normalization condition $\varphi(0)=1$, is equivalent to: \be
\left\{ \baa{rcl} \varphi''
& = &( k_z(\lambda) -m)\varphi  \hbox{ on }[0,l/2),\\
\varphi''
& = & k_z(\lambda)  \varphi  \hbox{ on }[l/2,l/2+z),\\
\varphi''
& = &( k_z(\lambda) -m)\varphi  \hbox{ on }[l/2+z,l+z),\\
\varphi''
& = & k_z(\lambda) \varphi  \hbox{ on }[l+z,L_0),\\
\varphi(0)& = & 1, \ \varphi(L_0)  =  e^{\lambda L_0} \varphi(0),
\
 \varphi'(L_0) = e^{\lambda L_0} \varphi'(0).
\eaa\right.\label{ax1}\ee

For each $z\in [0,L_0-l]$, let $\lambda^*_z$ be defined by the
formula (\ref{forc*_frag}). We have the following lemma:

\begin{lemm}
Assume that $l>3 L_0/4 $. Then, for all  $ \ds{z\in [0,L_0-l]}$,
we have $k_z(\lambda^*_z)>m$. \label{lem5_1}
\end{lemm}

\noindent\textbf{Proof of Lemma~\ref{lem5_1}.} Let us divide
equation (\ref{Lcmu}) by $\psi$ and integrate by parts over
$[0,L_0]$. Using the $L_0$-periodicity of $\psi$, we obtain:
$$\int_0^{L_0}\frac{|\psi'|^2}{\psi^2}+L_0 \lambda^2+\int_0^{L_0}
\mu_z(x)dx =L_0 k_z(\lambda).$$ Thus, \be k_z(\lambda)\geq
\lambda^2+ \frac{1}{L_0}\int_0^{L_0}\mu_z(x)dx=\lambda^2+ m
\frac{l}{L_0}. \label{eqlem51_1} \ee  From (\ref{forc*_frag}) and
(\ref{eqlem51_1})  we get:
$$(\lambda^*_z)^2+m \frac{l}{L_0}\leq k_z(\lambda^*_z)\leq
2\lambda^*_z\sqrt{m}.$$ Thus,
$(\lambda^*_z)^2-2\lambda^*_z\sqrt{m}+m l/L_0\leq 0$, which
 implies that $$\lambda^*_z\geq \sqrt{m}-\sqrt{m-m l/L_0}.$$Using
 (\ref{eqlem51_1}), we finally get $$k_z(\lambda^*_z)\geq 2m
 (1-\sqrt{1-l/L_0})>m,$$as soon as $l>3L_0/4$.\hfill{$\Box$}

\indent We now turn to the proof of Theorem~\ref{th_fragmentation}
and we assume that $l\in(3L_0/4,L_0)$. Using the fact that
$\varphi\in C^1(\R)$, a straightforward but lengthy computation
shows that, whenever $k_z(\lambda)>m$,  system (\ref{ax1}) is
equivalent to
$$\frac{F(z,\lambda,k_z(\lambda))}{G(z,k_z(\lambda))}=0,$$where $F$ and $G$ are two functions,
defined respectively  in  $[0,L_0-l]\times (0,+\infty) \times
[m,+\infty)$ and $[0,L_0-l]\times [m,+\infty)$ by: \be
\begin{array}{ll} F(z,\lambda,s) = &  4(2s-m) \sqrt
{s}\sqrt {s-m} \sinh(l\sqrt {s-m})\sinh(\alpha\sqrt {s})
\\ &
+ {m}^{2}\cosh(\beta\sqrt {s})(1-\cosh(l\sqrt{s-m}))
 \\ & +8(s^2-ms)[\cosh(l\sqrt {s-m})\cosh(\alpha\sqrt {s})-\cosh(\lambda L_0)]
\\&  +{m}^{2}\cosh(\alpha\sqrt
{s})(\cosh(l\sqrt{s-m})-1),
\end{array}  \label{FF}\ee
and \be
\begin{array}{ll} G(z,s) = & m\sqrt{s}\cosh(l\sqrt{s-m})\left[4 \sinh(\alpha\sqrt {s})(s/m-1) \right. \\ & \left.  +\left( \sinh(\alpha\sqrt {s}) - \sinh(\beta\sqrt {s}) \right)\left(1-1/\cosh(l\sqrt{s-m})\right)
 \right] \\ &+
m\sqrt{s-m}\sinh(l\sqrt{s-m})\cosh(\alpha
\sqrt{s})\left[\frac{4s}{m}-1+\frac{\cosh(\beta \sqrt
{s})}{\cosh(\alpha \sqrt{s})}\right],
\end{array} \label{G} \ee
with $\alpha:=L_0-l$ and $\beta:=L_0-l-2z$.

 Each factor in the
expression (\ref{G}) is positive, as soon as $s>m$, for $ \ds{z\in
[0,L_0-l]}$. Thus, whenever $k_z(\lambda)>m$, system (\ref{ax1})
is equivalent to the simpler equation \be
F(z,\lambda,k_z(\lambda))=0. \label{eqimplF} \ee Furthermore, from
Krein-Rutman theory, since the eigenfunction $\psi$ in
(\ref{Lcmu}) is positive, $k_z(\lambda)$ is the largest real
eigenvalue of the operator $\psi \mapsto \psi''+2\lambda\ \psi'
+\lambda^{2}\psi+\mu_z(x)\psi $. This result, implies that, for
each $ \ds{z\in [0,L_0-l]}$, and each $\lambda>0$,  $k_z(\lambda)$
is the largest real root of equation (\ref{eqimplF}), as soon as
$k_z(\lambda)>m$.

From equation (\ref{FF}), we easily see that \be \lim_{s\to
+\infty} F(z,\lambda,s)=+\infty, \label{lim+inf}\ee for all $
\ds{z\in [0,L_0-l]}$ and $\lambda>0$. Moreover, differentiating
(\ref{FF}) with respect to $z$, we obtain
$$\frac{\partial F}{\partial z}(z,\lambda,s)= 2 m^2 \sqrt{s}
\sinh(\sqrt{s}(L_0-l-2z))\left[\cosh(l\sqrt{s-m})-1\right].$$
Thus, for all $s>m$, and $\lambda>0,$ \be  \frac{\partial
F}{\partial z}(z,\lambda,s)>0 \hbox{ for }z\in [0,(L_0-l)/2),
\label{DFz2}\ee and $$\frac{\partial F}{\partial z}(z,\lambda,s)<0
\hbox{ for }z\in ((L_0-l)/2,L_0-l].$$

Now, take $z_1<z_2$ in $[0,(L_0-l)/2]$, and assume that
$c^*_{z_1}\leq c^*_{z_2}$.  It follows from formula
(\ref{forc*_frag}) that $k_{z_2}(\lambda)\geq c^*_{z_2}\lambda$,
for all $\lambda>0$. In particular, \be
k_{z_2}(\lambda^*_{z_1})\geq c^*_{z_2} \lambda^*_{z_1} \geq
c^*_{z_1} \lambda^*_{z_1} =k_{z_1}(\lambda^*_{z_1}).
\label{eqcrois}\ee From Lemma \ref{lem5_1}, we know that
$k_{z_1}(\lambda^*_{z_1})>m$. Thus, (\ref{eqcrois}) implies
$k_{z_2}(\lambda^*_{z_1})>m$. From the above discussion,
$k_{z_2}(\lambda^*_{z_1})$ is therefore the largest real root of
the equation $F(z_2,\lambda^*_{z_1},k_{z_2}(\lambda^*_{z_1}))=0$,
and, similarly, $k_{z_1}(\lambda^*_{z_1})$ is the largest real
root of $F(z_1,\lambda^*_{z_1},k_{z_1}(\lambda^*_{z_1}))=0$. Using
(\ref{lim+inf}) and (\ref{DFz2}), and since $0\le
z_1<z_2\le(L_0-l)/2$, we obtain
$k_{z_2}(\lambda^*_{z_1})<k_{z_1}(\lambda^*_{z_1})$, which
contradicts (\ref{eqcrois}). Therefore,  $c^*_z $ is
 a decreasing function of $z$ in $[0,(L_0-l)/2]$. Similar arguments
 imply that $c^*_z $ is an increasing function of $z$ in $[(L_0-l)/2,L_0-l]$. This concludes the proof of Theorem~\ref{th_fragmentation}.\hfill{$\Box$}

%%%%%%%%%%%%%%%%%%%%%%%%%%%%%%%%%%%%%%%%
%%%%%%%%%%%%%%%%%%%%%%%%%%%%%%%%%%%%%%%%

%%%%%%%%%%%%%%%%%%%%%%%%%%%%%%%%%%%%%%%%
%%%%%%%%%%%%%%%%%%%%%%%%%%%%%%%%%%%%%%%%

\end{document}